\documentclass[preprint]{elsarticle}

\usepackage{url}
\usepackage{color}
\usepackage{graphicx}
\usepackage{tikz,tikz-3dplot}

\usepackage{amsmath}
\usepackage{amsfonts}
\usepackage{amssymb}
\usepackage{amsthm}
\usepackage{epstopdf}
\usepackage{algorithmic}
\usepackage{setspace}
\usepackage[cm]{fullpage}
\usepackage{multirow}

\usepackage{fancyvrb}
\usepackage[colorinlistoftodos,prependcaption,textsize=tiny]{todonotes}
\usepackage{comment}

\usepackage[nomarkers]{endfloat}

\makeatletter  
\renewcommand{\verbatim@font}{%
  \ttfamily\small\catcode`\<=\active\catcode`\>=\active%
}
\makeatother   

\newcommand\restr[2]{{
  \left.\kern-\nulldelimiterspace 
  #1 
  \vphantom{\big|} 
  \right|_{#2} 
}}

\theoremstyle{remark}

\newtheorem{algorithm}{Algorithm}

\newcommand{\reviewalter}[1]{{\color{blue}#1}}

\def\stiff{\textrm{\textbf{K}}}

\newcommand{\mbf}[1]{\mbox{\boldmath$#1$}}

\journal{Applied Mathematics and Computation}

\begin{document}

\begin{frontmatter}

\title{Efficient and flexible MATLAB implementation of 2D and 3D elastoplastic problems}

\author[km,kam,enet,ugn]{M.~\v Cerm\'ak}
\ead{martin.cermak@vsb.cz}

\author[ugn]{S.~Sysala}
\ead{stanislav.sysala@ugn.cas.cz}

\author[bohemia,ascr]{J.~Valdman\corref{cor}}
\ead{jvaldman@prf.jcu.cz}

\cortext[cor]{Corresponding author}
\address[km]{Department of Mathematics, Faculty of Civil Engineering, V\v SB-TU Ostrava, Ostrava, Czech Republic}
\address[kam]{Department of Applied Mathematics, FEEIC, V\v SB-TU Ostrava, Ostrava, Czech Republic}
\address[enet]{ENET Centre, V\v SB-TU Ostrava, Ostrava, Czech Republic}
\address[ugn]{Institute of Geonics of the Czech Academy of Sciences, Ostrava, Czech Republic}
\address[bohemia]{Institute of Mathematics and Biomathematics, University of South Bohemia, \v Cesk\'e Bud\v ejovice, Czech Republic}
\address[ascr]{Institute of Information Theory and Automation of the Czech Academy of Sciences, Prague, Czech Republic}


\begin{abstract}
We propose an effective and flexible way to implement 2D and 3D elastoplastic problems in MATLAB using fully vectorized codes. Our technique is applied to a broad class of the problems including perfect plasticity or plasticity with hardening and several yield criteria. The problems are formulated in terms of displacements, discretized by the implicit Euler method in time and the finite element method in space, and solved by the semismooth Newton method. We discuss in detail selected models with the von Mises and Prager-Drucker yield criteria and four types of finite elements. The related codes are available for download. A particular interest is devoted to the assembling of tangential stiffness matrices. Since these matrices are repeatedly constructed in each Newton iteration and in each time step, we propose another vectorized assembling than current ones known for the elastic stiffness matrices. The main idea is based on a construction of two large and sparse matrices representing the strain-displacement and tangent operators, respectively, where the former matrix remains fixed and the latter one is updated only at some integration points. Comparisons with other available MATLAB codes show that our technique is also efficient for purely elastic problems. In elastoplasticity, the assembly times are linearly proportional to the number of integration points in a plastic phase and additional times due to plasticity never exceed assembly time of the elastic stiffness matrix. 
\end{abstract}


\end{frontmatter}


\section{Introduction}

The paper is focused on implementation of time-discretized elastoplastic problems formulated in terms of displacements. These problems include the following nonlinear variational equation \cite{NPO08, SH98} defined in each time step $k=1,2,\ldots, k_{max}$:
\begin{equation*}
(P_k)\quad \mbox{find } \mbf{u}_k\in\{\mbf u_{D,k}\}+\mathcal V: \quad \int_\Omega  T_k\left(\mbf\varepsilon(\mbf{u}_k)\right):\mbf\varepsilon(\mbf v)\, \mbox{d}x= \int_\Omega \mbf f_{V,k}.\mbf v\,\mbox{d}x+\int_{\Gamma_N} \mbf f_{t,k}.\mbf v\,\mbox{d}s \quad\forall \mbf v\in \mathcal V,
\label{eqn}
\end{equation*} 
where $\Omega\subset\mathbb R^d$, $d=2,3$, is a bounded domain with the Lipschitz boundary $\Gamma=\overline{\Gamma}_D\cup\overline{\Gamma}_N$. The parts $\Gamma_D$ and $\Gamma_N$ are open and disjoint. On $\Gamma_D$, the Dirichlet boundary conditions are prescribed and represented by a  given displacement function $\mbf u_{D,k}\in H^1(\Omega;\mathbb R^d)$. The testing functions $\mbf v$ belong to the space $\mathcal V$ of $H^1(\Omega;\mathbb R^d)$-functions vanishing on $\Gamma_D$. Further, 
$$\mbf f_{V,k}\in L^2(\Omega;\mathbb R^d), \qquad \mbf f_{t,k}\in L^2(\Gamma_N;\mathbb R^d)$$ are the prescribed volume and surface forces acting in $\Omega$ and on $\Gamma_N$, respectively. At any point 
$\mbf x\in\Omega:$ $$\mbf\varepsilon_k:=\mbf\varepsilon(\mbf{u}_k)=\frac{1}{2}(\nabla\mbf u_k+(\nabla \mbf u_k)^\top), \qquad \mbf\sigma_k:= T_k\left(\mbf\varepsilon(\mbf{u}_k)\right) $$ denote the infinitesimal strain tensor  $\mbf\varepsilon_k$ and the stress tensor $\mbf\sigma_k$  and $T_k:\mathbb R^{d\times d}_{sym}\rightarrow \mathbb R^{d\times d}_{sym}$ is a nonlinear stress-strain operator. The operator $T_k$ also depends on the plastic strain $\mbf{\varepsilon}^p_{k-1}$ and other internal variables known from the previous time step $k-1$. This operator varies depending on a particular constitutive model and is defined in an implicit form, in general. Although $T_k$ is nonsmooth, its semismoothness was proven for some elastoplastic models, see, e.g., \cite{GrVa, SaWi11, Sy14, SCL17}.


We follow a current computational procedure consisting of the following steps \cite{NPO08, SH98, SCKKZB16, SCL17}:
\begin{itemize}
\item[$(a)$] space discretization of $(P_k)$ by the finite element method (FEM); 
\item[$(b)$] solution of a resulting discretized system by the (semismooth) Newton method; 
\end{itemize}
Due to the possible presence of limit loads \cite{NPO08, HRS16a, HRS16b} leading to locking phenomena, implementations in elastoplasticity utilize higher order finite elements.  Their assemblies require suitable quadrature rules of higher order \cite{B06, NPO08}. 
To use the semismooth Newton method, one must find a generalized derivative of $T_k$ w.r.t. the strain variable, the so-called consistent tangent operator
 $$T^o_k:=T^o_k\left(\mbf\varepsilon(\mbf{u}_k)\right)=\mathrm{D}T_k\left(\mbf\varepsilon(\mbf{u}_k)\right).$$ After the solution $\mbf u_k$ is found, one can easily update $\mbf\sigma_k$, $\mbf{\varepsilon}^p_{k}$ and the internal variables at the level of integration points, and continue with the next time step.

 For a  computer implementation, it is crucial to find a suitable assembly of the so-called tangential stiffness matrix $\mbf K_{tangent}$  based on the operator $T^o_k$. We propose and explain in detail an efficient implementation of the tangential stiffness matrix in the form 
 \begin{equation}
\mbf K_{tangent}=\mbf K_{elast}+\mbf B^\top(\mbf D_{tangent}-\mbf D_{elast})\mbf B,
\label{K-split}
\end{equation}
 where  $\mbf{B}$ is a sparse matrix representing the strain-displacement operator at all integration points. The matrix $\mbf D_{tangent}$ is a block diagonal matrix and each block contains the operator $T^o_k$ for a particular integration point. The matrix $\mbf D_{elast}$ is also block diagonal and represents its elastic counterpart. It is applied to assembly the elastic stifness matrix  $\mbf K_{elast}$ in the form
 \begin{equation}
\mbf K_{elast}=\mbf B^\top\mbf D_{elast} \mbf B.
\label{K_elast-split}
\end{equation}
The matrices $\mbf K_{elast}, \mbf{B}, \mbf D_{elast}$ can be precomputed and only the matrix $\mbf D_{tangent}$ needs to be partially reassembled in each Newton iteration.  Although some MATLAB elastoplasticity codes are already available \cite{CK, Mcode, software_twoyield} and applied to various elastoplasticity models \cite{BCV05, SCKKZB16, SCL17}, they are not at all or only partially vectorized. The vectorization replaces time consuming loops by operations with long vectors and arrays and proves to be reasonably scalable and fast for large size problems. It typically takes only few minutes to solve studies elastoplastic benchmarks with several milion on unknows on current computers.  Authors are not aware of any other fully vectorized Matlab assembly of elastoplastic problems.

Our code is available for download \cite{software} and provides several computing benchmarks including
\begin{itemize}
\item[$(a)$] elastic and elastoplastic models with von Mises or Drucker-Prager  yield criteria; 
\item[$(b)$] finite element implementations of P1, P2, Q1, Q2 elements in both 2D and 3D.
\end{itemize}
Crucial functions are written uniformly regardless on these options. For the sake of brevity, we shall describe only 3D problems in this paper. Its plane strain reduction to 2D is usually straightforward and introduced within the code. Therefore, we set $d=3$ from now on. Although the paper is focused on an assembly of the tangential stiffness matrix, it is worth mentioning that the code contains complex implementation of elastic and elastoplastic problems defined on fixed geometries. In particular, there are included fully vectorized procedures for specific mesh generation, volume and surface forces. The solver is based on the semismooth Newton method combined with time stepping which can be adaptive, if necessary.

The rest of the paper is organized as follows. In Section \ref{sec_preliminaries}, we introduce a simplified scheme of elastoplastic constitutive problems and define the stress-strain operator $T_k$. Then we introduce examples of the operators $T_k$ and $T^o_k$ for some constitutive models. In Section \ref{sec_solution}, the finite element discretization of problem $(P_k)$ is described. In Section \ref{sec_algebra}, an algebraic formulation of the elastoplastic problem and the semismooth Newton method are introduced. In Section \ref{sec_notation}, basic MATLAB notation is introduced. In Sections \ref{sec_assembly} and \ref{sec_assembly_const}, we describe assembly of the elastic and tangent stiffness matrices and the vector of internal forces. In Section \ref{sec_example}, we illustrate the efficiency of the vectorized codes on particular 2D and 3D examples. The paper also contains Appendix, where reference elements, local basic functions and suitable quadrature formulas are summarized for the used finite elements.

\section{Elastoplastic constitutive model}
\label{sec_preliminaries}

In this section, we introduce an elastoplastic constitutive model and its implicit Euler discretization with respect to a time variable
$t\in[0, t_{\max}].$ The constitutive model is an essential part of the overall elastoplastic problem. It is defined at any point $\mbf x\in\Omega$. We shall assume that the model: satisfies the principle of maximum plastic dissipation; is based on linear elasticity; and contains optionally internal variables like kinematic or isotropic hardening. These assumptions enable to introduce variational formulation of the overall elastoplastic problem important for solvability analysis \cite{HR99}. More general constitutive models can be found, e.g.,  in \cite{NPO08, SH98}. 
The initial value constitutive problem has the following scheme \cite{CK, HR99, Sy14}:

\medskip\noindent
\textit{Given the history of the  infinitesimal strain tensor $\mbf\varepsilon=\mbf\varepsilon(t)\in\mathbb R^{3\times 3}_{sym}$, $t\in[0, t_{\max}]$, and the initial values 
$\mbf{\varepsilon}^p(0)=\mbf{\varepsilon}^p_0\in\mathbb R^{3\times 3}_{sym}, \; \mbf\xi(0)=\mbf\xi_0\in W$;
find the stress tensor $\mbf\sigma(t)\in\mathbb R^{3\times 3}_{sym}$, the plastic strain $\mbf{\varepsilon}^p(t)\in\mathbb R^{3\times 3}_{sym}$, internal variables $\mbf\chi(t)\in W$  and thermodynamical forces $\mbf \beta(t)\in W$ such that 
\begin{equation}
\left.
\begin{array}{c}
\mbf\sigma=\mathbb C(\mbf{\varepsilon}-\mbf{\varepsilon}^p),\;\; \mbf\beta=H(\mbf\chi),\\[1mm]
(\mbf\sigma,\mbf\beta)\in B,\\[1mm]
\dot{\mbf{\varepsilon}}^p:(\mbf\tau_\sigma-\mbf\sigma)+(-\dot{\mbf\beta},\mbf\tau_\beta-\mbf\chi)_W\leq0\quad\forall (\mbf\tau_\sigma,\mbf\tau_\beta)\in B
\end{array}
\right\}
\label{CIVP1}
\end{equation}
hold for each instant $t\in[0,t_{\max}]$, where $W$ is a finite dimensional space with the scalar product $(.,.)_W$, $H:W\rightarrow W$ is a Lipschitz continuous and strongly monotone function, $B\subset\mathbb R^{3\times3}_{sym}\times W$ is a closed convex set with a nonempty interior, and $\mathbb C: R^{3\times3}_{sym} \rightarrow R^{3\times3}_{sym}$ is the fourth order elastic tensor.}

\medskip\noindent
The elastic part of the strain tensor is denoted as
$$\mbf\varepsilon^e:=\mbf\varepsilon-\mbf\varepsilon^p$$  and the doubles $(\mbf\sigma,\mbf\beta)$, $(\mbf{\varepsilon}^p,-\mbf\beta)$ are called the generalized stress and strain, respectively. Notice that $(\dot{\mbf{\varepsilon}}^p,-\dot{\mbf\beta})$ belongs to the normal cone of $B$ at $(\mbf\sigma,\mbf\beta)$.  The definition of $W$ depends on the used internal variables. For example, one can set $W=\mathbb R^{3\times3}_{sym}$, $W=\mathbb R$, and $W=\{0\}$ for kinematic hardening, isotropic hardening, and perfect plasticity (no internal variables), respectively. 
Further, we shall assume that $\mathbb C$ is represented by two material parameters (isotropic material), e.g., by the bulk modulus $K>0$ and the shear modulus $G>0$ in the form
\begin{equation}
\mathbb C=K\mbf I\otimes\mbf I+2G\mathbb I_{D}.
\end{equation}
Here, $\mbf I\otimes\mbf I$ 
denotes the tensor product of unit (second order) tensors $\mbf I\in\mathbb R^{3\times 3}$ 
and $\mathbb I_{D}=\mathbb I-\frac{1}{3}\mbf I\otimes \mbf I$, where $\mathbb I\mbf\eta=\mbf\eta$ for any $\mbf\eta\in\mathbb R^{3\times 3}_{sym}$. Consequently, the Hooke's law states  the linear relation between the stress and the elastic strain tensors
\begin{equation}
\mbf\sigma=\mathbb C\mbf{\varepsilon}^e=K(\mbf I:\mbf{\varepsilon}^e)\mbf I+2G\mathbb I_{D}\mbf{\varepsilon}^e.
\label{elastic_law}
\end{equation}
The set $B$ is often in the form 
$$B=\{(\mbf\tau_\sigma,\mbf\tau_\beta)\in\mathbb R^{3\times3}_{sym}\times W\ |\; \Psi(\mbf\tau_\sigma,\mbf\tau_\beta)\leq 0\},$$
where $\Psi$ is a convex yield function. This fact allows to rewrite  (\ref{CIVP1}) by the Karush-Kuhn-Tucker (KKT) conditions:
\begin{equation}
\left.
\begin{array}{c}
\mbf\sigma=\mathbb C(\mbf{\varepsilon}-\mbf{\varepsilon}^p),\;\; \mbf\beta=H(\mbf\chi),\\[1mm]
(\dot{\mbf{\varepsilon}}^p,-\dot{\mbf\chi})\in\dot\lambda\partial \Psi(\mbf\sigma,\mbf\beta),\\[1mm]
\dot\lambda\geq0,\;\; \Psi(\mbf\sigma,\mbf\beta)\leq0,\;\; \dot\lambda \Psi(\mbf\sigma,\mbf\beta)=0.
\end{array}
\right\}
\label{CIVP2}
\end{equation}
Here, $\dot\lambda$ stands for the plastic multiplier and $\partial$ is a subdifferential operator. The formulation (\ref{CIVP2}) is broadly used in engineering practice \cite{NPO08, SH98} and it is convenient for finding analytical or semianalytical solution of the constitutive problem \cite{SCKKZB16, SCL17}.

\subsection{The implicit discretization of the constitutive problem}
\label{subsec_implicit}

Let us consider a partition of the time interval $[0,t_{\max}]$ in the form $$0=t_0<t_1<\ldots<t_k<\ldots<t_N=t_{\max}$$ and denote $\mbf{\sigma}_k:=\mbf{\sigma}(t_k)$, $\mbf\varepsilon_k:=\mbf\varepsilon(t_k)$, $\mbf\varepsilon^p_k:=\mbf\varepsilon^p(t_k)$, $\mbf\chi_k:=\mbf\chi(t_k)$, $\mbf\beta_k=\mbf\beta(t_k)$. The $k$-th step of the incremental constitutive problem discretized by the implicit Euler method reads as: 

\medskip\noindent
{\it Given $\mbf{\varepsilon}_k$, $\mbf\varepsilon^p_{k-1}$ and $\mbf\chi_{k-1}$, find $\mbf{\sigma}_k$,  $\mbf\varepsilon^p_k$, $\mbf\chi_k$, $\mbf\beta_k$, and $\triangle\lambda$ satisfying:}
\begin{equation}
\left.
\begin{array}{c}
\mbf\sigma_k=\mathbb C(\mbf{\varepsilon}_k-\mbf{\varepsilon}^p_k),\;\; \mbf\beta_k=H(\mbf\chi_k),\\[1mm]
(\mbf{\varepsilon}^p_k-\mbf{\varepsilon}^p_{k-1},-\mbf\chi_k+\mbf\chi_{k-1})\in\triangle\lambda\partial \Psi(\mbf\sigma_k,\mbf\beta_k),\\[1mm]
\triangle\lambda\geq0,\;\; \Psi(\mbf\sigma_k,\mbf\beta_k)\leq0,\;\; \triangle\lambda\Psi(\mbf\sigma_k,\mbf\beta_k)=0.
\end{array}
\right\}
\label{CIVP_k}
\end{equation}
Under the assumptions mentioned above, problem (\ref{CIVP_k}) has a unique solution \cite{Sy14}. Therefore, one can define the stress-strain operator $T_k$ as follows:
$$\mbf\sigma_k=T_k(\mbf\varepsilon_k)=T(\mbf\varepsilon_k;\mbf\varepsilon^p_{k-1},\mbf\chi_{k-1}).$$ 
In general, the function $T_k$ is implicit. Nevertheless, it is well-known that problem (\ref{CIVP_k}) can be simplified and, sometimes, its solution can be found in a closed form. To this end the elastic predictor -- plastic corrector method is applied. Within the elastic prediction, it is checked whether the trial generalized stress $(\mbf\sigma^{tr}_k,\mbf\beta_{k-1})$, $\mbf\sigma^{tr}_k=\mathbb C(\mbf\varepsilon_k-\mbf\varepsilon^p_{k-1})$ is admissible or not. If $\Psi(\mbf\sigma^{tr}_k,\mbf\beta_{k-1})\leq0$ then 
$$\triangle\lambda=0,\;\;\mbf\sigma_k=\mbf\sigma_k^{tr}, \;\; \mbf\varepsilon^p_{k}=\mbf\varepsilon^p_{k-1},\;\;\mbf\chi_{k}=\mbf\chi_{k-1},\;\; \mbf\beta_{k}=\mbf\beta_{k-1}$$
is the solution to (\ref{CIVP_k}) representing the elastic response. Otherwise, (\ref{CIVP_k})$_3$ reduces into $\triangle\lambda>0$ and $\Psi(\mbf\sigma_k,\mbf\beta_k)=0$, and the plastic correction (return mapping) of the trial generalized stress is necessary.
From \cite{Sy14}, it follows that $T_k$ is Lipschitz continuous in $\mathbb R^{3\times 3}_{sym}$. Therefore, one can define a function 
$$T^o_k\colon\mathbb R^{3\times 3}_{sym}\rightarrow \mathcal L({\mathbb R^{3\times 3}_{sym},\mathbb R^{3\times 3}_{sym}})$$ representing a generalized Clark derivative of $T_k$. Clearly, if $T_k$ is differentiable at $\mbf\varepsilon_k$ then $T_k^o(\mbf\varepsilon_k)=\mathrm{D}T_k(\mbf\varepsilon_k)$. One can also investigate the semismoothness of $T_k$ under the assumption that $H$ is semismooth in $W$ \cite{Sy14}. 

In the elastic case, no plastic strain and internal variables occur, $\mbf\varepsilon^p_k=0, \chi_k=0$, the stress-strain operator and its derivative simplify as 
\begin{equation}
T_k(\mbf\varepsilon_k)=\mathbb C\mbf\varepsilon_k,\quad T_k^o(\mbf\varepsilon_k)=\mathbb C.
\label{elastic_operator}
\end{equation}

Further, we introduce two particular examples of elastoplastic models and their operators $T_k$ and $T_k^o$, for illustration.

\begin{figure}
\center
\includegraphics[width=0.3\textwidth]{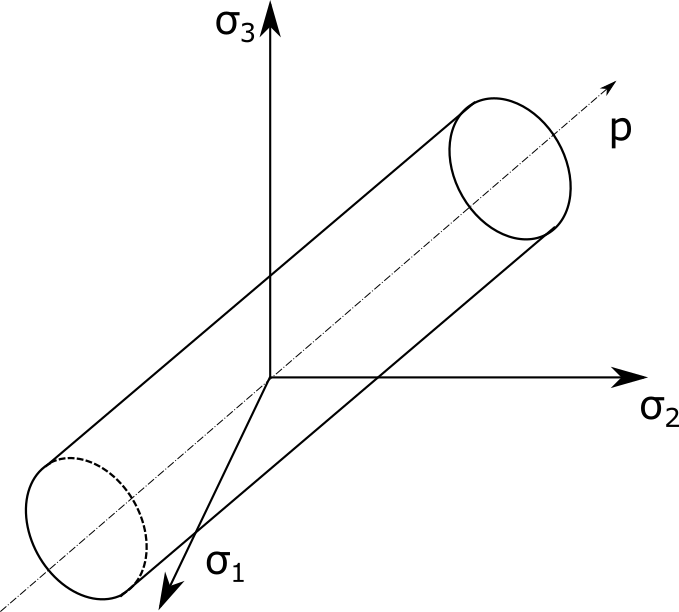}
\hspace*{10mm}
\includegraphics[width=0.3\textwidth]{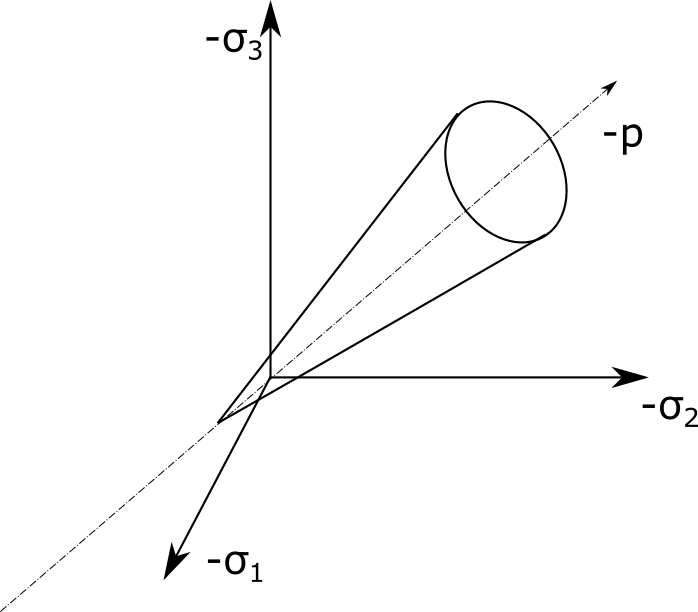}
\caption{The von Mises yield criterion (left) and the Drucker-Prager yield criterion (right). $\sigma_1$, $\sigma_2$, $\sigma_3$ denote the principal stresses and $p=(\sigma_1+\sigma_2+\sigma_3)/3=\mbf I:\mbf\sigma/3$ is the hydrostatic pressure.}
\label{fig_yield_criteria}

\end{figure}
\subsection{Von Mises yield criterion and linear kinematic hardening}
\label{subsec_kinematic}

This model corresponds to the choice 
$$W=\mathbb R^{3\times3}_{sym}, \quad \mbf\beta=H(\mbf \chi):=a\mbf\chi, \quad a>0, \quad \Psi(\mbf\tau_\sigma,\mbf\tau_\beta)=|\mathbb I_D(\mbf\tau_\sigma-\mbf\tau_\beta)|- Y, \quad Y>0,$$ where  $\mathbb I_{D}\mbf\tau$ is the deviatoric part of $\mbf\tau\in \mathbb R^{3\times3}_{sym}$. If $\beta=0$ then the set $B$ of  admissible stress tensors is a cylinder aligned with the hydrostatic axis, see Figure \ref{fig_yield_criteria}. The kinematic hardening causes translation of the cylinder in the normal direction to the yield surface.
The corresponding solution of (\ref{CIVP_k}) can be found, e.g., in \cite{Bl99, CK}. We arrive at: 
\begin{equation}
T_k(\mbf\varepsilon_k)=\left\{
\begin{array}{ll}
\mbf\sigma_k^{tr},& |\mbf s_k^{tr}|\leq Y,\\
\mbf\sigma_k^{tr}-\frac{2G}{2G+a}\left(|\mbf s_k^{tr}|- Y\right)\mbf n_k^{tr},& |\mbf s_k^{tr}|> Y,
\end{array}
\right.
\label{T_kinematic}
\end{equation}
\begin{equation}
T_k^o(\mbf\varepsilon_k)=\left\{
\begin{array}{ll}
\mathbb C,& |\mbf s_k^{tr}|\leq Y,\\
\mathbb C-\frac{4G^2}{2G+a}\mathbb I_D+\frac{4G^2}{2G+a}\frac{ Y}{|\mbf s_k^{tr}|}\left(\mathbb I_D-\mbf n_k^{tr}\otimes\mbf n_k^{tr}\right),& |\mbf s_k^{tr}|> Y,
\end{array}
\right.
\label{DT_kinematic}
\end{equation}
where $$\mbf\sigma^{tr}_k=\mathbb C(\mbf\varepsilon_k-\mbf\varepsilon^p_{k-1}),\;\; \mbf s_k^{tr}=\mathbb I_D\mbf\sigma_k^{tr}-\mbf\beta_{k-1},\;\;\mbf n_k^{tr}=\frac{\mbf s_k^{tr}}{|\mbf s_k^{tr}|}.$$
Further, the hardening variable and the plastic strain are updated as follows:
\begin{equation}
\mbf\beta_k=\left\{
\begin{array}{ll}
\mbf\beta_{k-1},& |\mbf s_k^{tr}|\leq Y,\\
\mbf\beta_{k-1}+\frac{a}{2G+a}\left(|\mbf s_k^{tr}|- Y\right)\mbf n_k^{tr},& |\mbf s_k^{tr}|> Y,
\end{array}
\right.\quad
\label{T_other}
\mbf\varepsilon^p_k=\left\{
\begin{array}{ll}
\mbf\varepsilon^p_{k-1},& |\mbf s_k^{tr}|\leq Y,\\
\mbf\varepsilon^p_{k-1}+\frac{1}{2G+a}\left(|\mbf s_k^{tr}|- Y\right)\mbf n_k^{tr},& |\mbf s_k^{tr}|> Y.
\end{array}
\right.
\end{equation}

Notice that if $a=0$ we arrive at formulas for the perfect plastic model with the von Mises yield criterion.

\subsection{Drucker-Prager yield criterion and perfect plasticity}
\label{subsec_DP}

In perfect plasticity, the internal variables $\mbf\beta$ and $\mbf\chi$ are not included in the model and $W=\{0\}$. The Drucker-Prager yield function is defined as follows:
$$\Psi(\mbf\tau_\sigma)=\sqrt{\frac{1}{2}}|\mathbb I_D\mbf\tau_\sigma|+\frac{\eta}{3} \mbf I:\mbf\tau_\sigma-c,\quad \mbf\tau\in\mathbb R^{3\times 3}_{sym}, $$
where $\eta, c>0$ are given material parameters. The corresponding set $B$ of admissible stress tensors is depicted in Figure \ref{fig_yield_criteria}. It is a cone aligned with the hydrostatic axis. The operators $T_k$ and $T_k^o$ can be found in closed forms \cite{NPO08, SCKKZB16} by solving the system (\ref{CIVP_k}). To summarize their forms, we distinguish three different cases and define the following auxiliary notation:
$$\mbf\sigma^{tr}_k=\mathbb C(\mbf\varepsilon_k-\mbf\varepsilon^p_{k-1}),\;\; p_k^{tr}=\mbf I:\mbf\sigma^{tr}_k,\;\; \mbf s_k^{tr}=\mathbb I_D\mbf\sigma_k^{tr},\;\; \varrho_k^{tr}=|\mbf s_k^{tr}|,\;\;\mbf n_k^{tr}=\frac{\mbf s_k^{tr}}{|\mbf s_k^{tr}|}.$$

\begin{enumerate}
\item {\it Elastic response} occurs if $\Psi(\mbf\sigma^{tr}_k)\leq0$. Then $T_k(\mbf\varepsilon_k)=\mbf\sigma^{tr}_k$, $\mbf\varepsilon^p_k=\mbf\varepsilon^p_{k-1}$ and $T_k^o(\mbf\varepsilon_k)=\mathbb C$.
\item {\it Return to the smooth portion of the yield surface} occurs if $\Psi(\mbf\sigma^{tr}_k)>0$ and $\eta p_k^{tr}-\frac{K\eta^2}{G\sqrt{2}}\varrho_k^{tr}<c$. Then,
\begin{equation*}
T_k(\mbf\varepsilon_k)=\mbf{\sigma}^{tr}_k-\frac{\Psi(\mbf\sigma^{tr}_k)}{G+K\eta^2}\left(G\sqrt{2}\mbf{n}^{tr}_k+K\eta\mbf I\right),\quad 
\mbf\varepsilon^p_k=\mbf\varepsilon^p_{k-1}+\frac{\Psi(\mbf\sigma^{tr}_k)}{G+K\eta^2}\left(\frac{\sqrt{2}}{2}\mbf{n}^{tr}_k+\frac{1}{3}\eta\mbf I\right)
\end{equation*}
\begin{equation*}
T_k^o(\mbf\varepsilon_k)=\mathbb C-\frac{1}{G+K\eta^2}\left(\frac{2G^2\sqrt{2}\Psi(\mbf\sigma^{tr}_k)}{\varrho^{tr}}\left(\mathbb I_{D}-\mbf{n}^{tr}_k\otimes\mbf{n}^{tr}_k\right)+(G\sqrt{2}\mbf{n}^{tr}_k+K\eta\mbf I)\otimes (G\sqrt{2}\mbf{n}^{tr}_k+ K\eta\mbf I)\right).
\end{equation*}
\item {\it Return to the apex of the yield surface} occurs if $\eta p_k^{tr}-\frac{K\eta^2}{G\sqrt{2}}\varrho_k^{tr}\geq c$. Then,
\begin{equation*}
T_k(\mbf\varepsilon_k)=\frac{c}{\eta}\mbf I,\quad T_k^o(\mbf\varepsilon_k)=\mathbb O,\quad \mbf\varepsilon^p_k=\mbf\varepsilon_k-\frac{c}{3K\eta}\mbf{I},
\end{equation*}
where $\mathbb O$ is a zero fourth order tensor, i.e., $\mathbb O\mbf\tau:\mbf\zeta=0$ for any $\mbf\tau,\mbf\zeta\in\mathbb R^{3\times 3}_{sym}$.
\end{enumerate}

\section{Finite element discretization}
\label{sec_solution}

The standard Galerkin method leads to the following discrete counterpart of problem $(P_k)$:
\begin{equation*}
(P_k)_h\quad \mbox{find } \mbf{u}_{k,h}\in\{\mbf u_{D,k,h}\}+\mathcal V_h: \quad \int_{\Omega_h}  T_k\left(\mbf\varepsilon(\mbf{u}_{k,h})\right):\mbf\varepsilon(\mbf v_h)\, \mbox{d}x= \int_{\Omega_h} \mbf f_{V,k,h}.\mbf v_h\,\mbox{d}x+\int_{\Gamma_{N,h}} \mbf f_{t,k,h}.\mbf v_h\,\mbox{d}s \quad\forall \mbf v_h\in \mathcal V_h,
\label{eqn_h}
\end{equation*}
where $\mathcal V_h$ is a finite dimensional approximation of $\mathcal V$ and $\mbf f_{V,k,h}$, $\mbf f_{t,k,h}$, $\mbf u_{D,k,h}$, and $\Omega_h$ are suitable approximations of $\mbf f_{V,k}$, $\mbf f_{t,k}$, $\mbf u_{D,k}$, and $\Omega$, respectively. Beside these approximations, one must also take into account a numerical integration specified below. The finite element method corresponds to a specific choice of $\mathcal V_h$. In particular, we shall consider conforming and isoparametric finite elements of Lagrange type in 3D and define the corresponding basis functions of $\mathcal V_h$. For more detail, we refer to \cite{B06}. 

Let $\hat T\subset\mathbb R^3$ denote a reference element  w.r.t. to Cartesian coordinates $\mbf\xi=(\xi_1,\xi_2,\xi_3)$. The reference element is usually a convex polyhedron (in our case, either a tetrahedron or a hexahedron) and its boundary consists of sides $\hat E\subset\partial \hat T$. On $\hat T$, we prescribe: a space $\hat S$ of dimension $n_p$, basis functions $\hat\Phi_p=\hat\Phi_p(\mbf\xi)$ of $\hat S$, and nodes $\hat N_p\in \hat T$, $p=1,2,\ldots,n_p$ such that $\hat\Phi_p(\hat N_q)=\delta_{pq}$, $p,q=1,2,\ldots,n_q$. Further,  the domain $\bar\Omega_h$ is covered by a regular triangulation $\mathcal T_h$, i.e., $\bar\Omega_h=\bigcup_{T\in\mathcal T_h}T$. We assume that any element $T$ can be described by nodes $N_{T,p}\in T$, $p=1,2,\ldots,n_p$, and by the following nondegenerative transformation of $\hat T$:
\begin{equation}
\forall\mbf x\in T\ \exists! \mbf\xi\in\hat T:\quad \mbf x=\sum_{p=1}^{n_p} \hat \Phi_p(\mbf\xi) N_{T,p}.
\label{transformation}
\end{equation}
Then the so-called iso-parametric transformation (\ref{transformation}) maps  $\hat N_{T,p}$ onto $N_{T,p}$ for $p=1,2,\ldots, n_p$ and any side $\hat E\subset\partial \hat T$ has a corresponding side $E\subset\partial T$. It is required that: any side $E$ is either an intersection of two neighboring elements or a part of the boundary $\partial \Omega_h$; if $N_{T,p}\in E$, $E=T\cap T'$, then $N_{T,p}$ is also a nodal point of $T'\in \mathcal T_h$. Further, we standardly assume that if some node $N_{T,p}$ belongs to $\partial \Omega_h$ then also $N_{T,p}\in\partial \Omega$ in order to reflect the curvature of the original domain $\Omega$. 

Making use of the transformation  (\ref{transformation}), one can define the local basis function $\Phi_{T,p}$, $p=1,2,\ldots, n_p$ for any element $T$:
\begin{equation}
\Phi_{T,p}(\mbf x):=\hat\Phi_{p}(\mbf \xi),\quad\mbf x\in T.
\label{local_base}
\end{equation}
Define the set $\mathcal N$ of all nodes corresponding to the triangulation $\mathcal T_h$, i.e., $N_{T,p}\in\mathcal N$ for any $T\in\mathcal T_h$ and $p=1,2,\ldots,n_p$. Let $n_n$ denote a number of nodes from $\mathcal N$. For any node $N_j\in\mathcal N$, $j=1,2,\ldots, n_n$, we define the function $\Phi_{j}\colon\Omega_h\rightarrow \mathbb R$  as follows:
\begin{equation}
\Phi_{j}|_T:=\left\{
\begin{array}{cc}
\Phi_{T,p},& \mbox{if }\exists p\in\{1,2,\ldots,n_p\}:\; N_j=N_{T,p},\\[2mm]
0, & \mbox{otherwise},
\end{array}\right.\quad\forall T\in\mathcal T_h.
\label{base_function}
\end{equation}
The assumption on conforming finite elements means that the space $\hat S$ and the reference nodes $\hat N_p\in \hat T$, $p=1,2,\ldots,n_p$, are such that the functions $\Phi_{j}$, $j=1,2,\ldots, n_n$, are continuous. This holds, e.g., for simplicial P1, P2, $\ldots$ or quadrilateral Q1, Q2, $\ldots$ elements, $\ell=1,2,\ldots$, see \cite{B06}. Let $\mathcal H_h$ denote the space generated by (global) basis functions $\Phi_{j}$, $j=1,2,\ldots, n_n$, i.e., $\mathcal H_h=\mathrm{lin}\{\Phi_1, \Phi_2,\ldots, \Phi_{n_n}\}$. It is a space of continuous and piecewise smooth functions of the form
\begin{equation}
v_h(\mbf x)=\sum_{j=1}^{n_n} \Phi_j(\mbf x)v_h(N_j),\quad \mbf x\in\Omega_h.
\label{test_function}
\end{equation}
The space $\mathcal H_h$ is a finite dimensional subspace of $H^1(\Omega)$ and analogously, $\mathcal H_h^3:=\mathcal H_h\times \mathcal H_h\times\mathcal H_h\subset H^1(\Omega;\mathbb R^3)$. The required space $\mathcal V_h$ is the subspace of $\mathcal H_h^3$ which does not contain basis functions corresponding to the nodes lying on $\Gamma_D$, i.e., 
$$\mathcal V_h=\mathrm{lin}\{(\Phi_{j},0,0), (0, \Phi_{j},0), (0,0,\Phi_{j});\; j\in \mathcal I_Q\},$$ where $\mathcal I_Q=\{j\in\{1,2,\ldots,n_n\}\ |\; N_{j}\not\in\Gamma_D\}$. We shall also use the notation $\mathcal I_D:=\{1,2,\ldots,n_n\}\setminus\mathcal I_Q$.

In order to evaluate volume integrals, we standardly split the domain $\Omega_h$ into elements $T\in\mathcal T_h$ and use the transformation (\ref{transformation}) of $T$ onto $\hat T$. We have the following transforming formulas \cite{B06}:
\begin{eqnarray}
\mathrm d\mbf x=|\det J_T(\mbf\xi)|\,\mathrm d\mbf\xi, \quad \left(\begin{array}{l}\frac{\partial}{\partial x_1}\\ \frac{\partial}{\partial x_2}\\ \frac{\partial}{\partial x_3}\end{array}\right)= J_T(\mbf\xi)^{-1}\left(\begin{array}{l}\frac{\partial}{\partial \xi_1}\\ \frac{\partial}{\partial \xi_2}\\ \frac{\partial}{\partial \xi_3}\end{array}\right),\label{Jacobian2}\\
\quad [J_T(\mbf\xi)]_{i,j}=\frac{\partial x_j}{\partial \xi_i}\stackrel{ (\ref{transformation})}{=}\sum_{p=1}^{n_p} \frac{\hat \Phi_p(\mbf\xi)}{\partial \xi_i} [N_{T,p}]_j,\quad i,j=1,2,3,
\label{Jacobian}
\end{eqnarray}
where $J_T$ denotes a Jacobian matrix (referred to as Jacobian) on $T\in\mathcal T_h$. Finally, we consider a numerical quadrature in $\hat T$ and write it in an abstract form:
\begin{equation}
\int_{\hat T} g(\mbf\xi)\,\mathrm d\mbf\xi \approx \sum_{q=1}^{n_q}\omega_q g(\hat A_q),
\label{quadrature}
\end{equation}
where $\hat A_q\in\hat T$ are the quadrature points and $\omega_q$ are the corresponding weights for $q=1,2,\ldots,n_q$. For example, we shall use the following formulas:
\begin{equation}
\int_{T} f\Phi_{T,p}\,\mathrm d\mbf x \approx \sum_{q=1}^{n_q}\omega_q
|\det J_T(\hat A_q)| f(A_{T,q})\hat\Phi_p(\hat A_q),\quad p=1,2,\ldots,n_p,
\label{formula1}
\end{equation}
\begin{equation}
\int_{T} f\frac{\partial\Phi_{T,p}}{\partial x_i}\,\mathrm d\mbf x \approx \sum_{q=1}^{n_q}\omega_q
|\det J_T(\hat A_q)| f(A_{T,q})\frac{\partial \Phi_{T,p}(A_{T,q})}{\partial x_i}, \quad p=1,2,\ldots,n_p,\; i=1,2,3,
\label{formula2}
\end{equation}
where the nodes $A_{T,q}\in T$ correspond to $\hat A_q$ within the transformation (\ref{transformation}) and 
\begin{equation}
\frac{\partial \Phi_{T,p}(A_{T,q})}{\partial x_i}\stackrel{(\ref{Jacobian2})}{=}\sum_{j=1}^3[J_T(\hat A_q)^{-1}]_{i,j}\frac{\partial\hat\Phi_p(\hat A_q)}{\partial \xi_j}, \qquad p=1,2,\ldots,n_p,\; i=1,2,3.
\label{basis_derivative}
\end{equation}
From these formulas, one can easily derive the required volume integrals through $\Omega_h$ important for assembly of the stiffness matrix and of the vectors of volume forces (internal and external), see Sections \ref{sec_assembly} and \ref{sec_assembly_const}. Examples of finite elements and convenient numerical quadratures are introduced in the Appendix.

\section{Algebraic problem and Newton-like method}
\label{sec_algebra}

An algebraic form of problem $(P_k)_h$ reads as:
\begin{equation*}
(\mathbb P_k)\qquad \mbox{find } \mathbf{u}_{k}\in\{\mathbf u_{D,k}\}+\mathbb V: \quad \mathbf v^\top(F_k(\mathbf u_k)-\mathbf f_k)=0  \quad\forall \mathbf v\in \mathbb V.
\label{eqn_alg}
\end{equation*}
Here, $\mathbf v$, $\mathbf u_k$ and $\mathbf u_{D,k}$ are the algebraic counterparts of $\mbf v_h$, $\mbf u_{k,h}$ and $\mbf u_{D,k,h}$, respectively. For example,
$$\mathbf v=(v_{h,1}(N_1), v_{h,2}(N_1), v_{h,3}(N_1),\ldots, v_{h,1}(N_{n_n}), v_{h,2}(N_{n_n}), v_{h,3}(N_{n_n}))^\top\in\mathbb R^{3n_n}.$$
The remaining notation is defined as follows:
$$\mathbb V:=\{\mathbf v\in\mathbb R^{3n_n}\ |\; v_{3j-2}=v_{3j-1}=v_{3j}=0,\;j\in\mathcal I_D\},$$
$$\mathbf f_k\in\mathbb R^{3n_n},\quad \mathbf v^\top\mathbf f_k:= \int_{\Omega_h} \mbf f_{V,k,h}.\mbf v_h\,\mbox{d}x+\int_{\Gamma_{N,h}} \mbf f_{t,k,h}.\mbf v_h\,\mbox{d}s\quad\forall \mbf v_h\in\mathcal V_h,$$
$$F_k\colon\mathbb R^{3n_n}\rightarrow\mathbb R^{3n_n},\quad\mathbf v^\top F_k(\mathbf u_k):=\int_{\Omega_h}  T_k\left(\mbf\varepsilon(\mbf{u}_{k,h})\right):\mbf\varepsilon(\mbf v_h)\, \mbox{d}x\quad\forall \mbf v_h\in\mathcal V_h.$$
Problem $(\mathbb P_k)$ can be simply transformed to a system of nonlinear equations by elimination of rows that correspond to the Dirichlet boundary conditions:
\begin{equation*}
\mbox{find } \mathbf{u}_{k}\in\{\mathbf u_{D,k}\}+\mathbb V: \quad \tilde F_k(\mathbf u_k)=\tilde{\mathbf f}_k\quad\mbox{in }\mathbb R^n,\;n=\mathrm{dim}\,\mathbb V,
\label{eqn_alg2}
\end{equation*}
where $\tilde F_k$ and $\tilde{\mathbf f}_k$ arise from $F_k$ and $\mathbf f_k$, respectively, by the elimination.  

In order to introduce the semismooth Newton method, we define a generalized gradient of $F_k$ at $\mathbf u$:
$$K_k\colon\mathbb R^{3n_n}\rightarrow\mathbb R^{3n_n\times 3n_n},\quad\mathbf v^\top K_k(\mathbf u)\mathbf w:=\int_{\Omega_h}  T_k^o\left(\mbf\varepsilon(\mbf{u}_{k,h})\right)\mbf\varepsilon(\mbf w_h):\mbf\varepsilon(\mbf v_h)\, \mbox{d}x\quad\forall \mbf v_h, \mbf w_h\in\mathcal V_h.$$
The corresponding algorithm reads as:


\begin{algorithm}[ALG-NEWTON]
\hspace{0.2cm}
\begin{spacing}{1.2}
\begin{algorithmic}[1]
  \STATE initialization: $\mathbf{u}_{k}^0=\mathbf u_{D,k}$
  \FOR{$\ell=1,2,\ldots$}
    \STATE find $\delta \mathbf{u}^{\ell}\in\mathbb V$: $\;\mathbf v^\top\mathbf {K}_k^\ell\delta  \mathbf{u}^{\ell}=\mathbf v^\top[\mathbf{f}_k-\mathbf {F}_{k}^\ell]$, $\;\;\forall \mathbf v\in \mathbb V$
    \STATE compute $\mathbf{u}_{k}^{\ell}=\mathbf{u}_{k}^{\ell-1}+\delta \mathbf{u}^{\ell}$
    \STATE {\bf if
    }{$\|\delta
    \mbf{u}^{\ell}\|_e/(\|\mathbf{u}_{k}^{\ell-1}\|_e+\|\mathbf{u}_{k}^\ell\|_e)\leq\epsilon_{\mathrm{Newton}}$} {\bf then stop}
  \ENDFOR
  \STATE set $\mathbf{u}_{k}=\mathbf{u}_{k}^{\ell}$,
\end{algorithmic}
\end{spacing}
\end{algorithm}

\noindent
where $\mathbf {K}_k^\ell:=K_k(\mathbf u_k^\ell)$, $\mathbf {F}_k^\ell:=F_k(\mathbf u_k^\ell)$ and $\|\mathbf v\|_e^2:=\mathbf v^\top\mathbf K_{elast}\mathbf v$ for any $\mathbf v\in\mathbb V$. In each $\ell$-th Newton iteration, we solve the linear problem with the tangent stiffness matrix $\mathbf {K}_k^\ell$. This problem can be transformed to the linear system of equations
\begin{equation}
\tilde{\mathbf K}_k^\ell\delta\tilde{\mathbf u}^\ell=\tilde{\mathbf f}_k-\tilde{\mathbf F}_k^\ell\quad\mbox{in }\mathbb R^n,\;n=\mathrm{dim}\,\mathbb V,
\label{lin_syst}
\end{equation}
by the elimination of the Dirichlet nodes. Assembly of $K_k(\mathbf u)$ and $F_k(\mathbf u)$ for some $\mathbf u\in\mathbb R^{3n_n}$ is derived in detail in Sections \ref{sec_assembly} and \ref{sec_assembly_const} where the indices $k$ and $\ell$ will be omitted, for the sake of simplicity.

It is well known that convergence of the Newton method is superlinear under the assumption that $\mathbf{u}_{k}^0=\mathbf u_{D,k}$ is close to the solution of $(\mathbb P_k)$. Alternatively, one can apply a damped version of the Newton method in elastoplasticity, see, e.g., \cite{Sy12}.


\section{Basic MATLAB notation}
\label{sec_notation}

 To distinguish notation for MATLAB commands and expressions in the next sections, a typescript is used (e.g., \verb+ELEM, COORD+, etc.). Instead of subscripts, we use the underscore symbol in MATLAB, e.g., \verb+n_n, n_e, n_q+, etc. Further, to be the codes vectorized (without long for-cycles), we work with arrays and use standard MATLAB commands like 
\begin{verbatim}
 repmat, kron, reshape, .*, ./, sum.
\end{verbatim}

Let $n_n$ and $n_e$ be numbers of nodes from $\mathcal N$ and elements from $\mathcal T_h$, respectively, and recall that $n_p$ and $n_q$ denote number of nodes and quadrature points within a finite element, respectively. Further, we define $n_{int}:=n_en_q$, i.e., a number of all integration points $A_{T,q}$, $T\in\mathcal T_h$, $q\in\{1,\ldots,n_q\}$. 

We use a $3\times n_n$ array \verb+U+ for storage of nodal displacements. In order to receive the corresponding displacement vector $\mathbf u\in\mathbb R^{3n_n}$, it suffices to use the command \verb+U(:)+, i.e., $\mathbf u=$\verb+U(:)+. Further, it is convenient to define the logical $3\times n_n$ array \verb+Q+ which indicates the nodes belonging to $\mathcal I_Q$, i.e., the nodes where the Dirichlet boundary condition is not prescribed. Then the restricted displacement vector $\tilde{\mathbf u}$ satisfies $\tilde{\mathbf u}=$\verb+U(Q)+. Similarly, one can restrict the stiffness matrix $\mbf K$ and the vectors of external and internal forces, $\mbf f$ and $\mbf F$. Once these objects are at disposal in MATLAB, the system (\ref{lin_syst}) of linear equations can be solved by the following commands (omitting indices $k$ and $\ell$):
\begin{verbatim}
 dU=zeros(3,n_n); dU(Q) = K(Q,Q)\(f(Q)-F(Q));
\end{verbatim}

Further, stress and strain tensors, $\mbf\sigma\in\mathbb R^{3\times 3}_{sym}$ and $\mbf \varepsilon\in\mathbb R^{3\times 3}_{sym}$, are represented standardly by the following vectors:
$$(\sigma_{11}, \sigma_{22}, \sigma_{33}, \sigma_{12}, \sigma_{23}, \sigma_{31})^\top,\quad (\varepsilon_{11}, \varepsilon_{22}, \varepsilon_{33}, 2\varepsilon_{12}, 2\varepsilon_{23}, 2\varepsilon_{31})^\top,$$ 
respectively\footnote{Within implementation, it is necessary to keep the different representations of stress-based and strain-based tensors.}.
These vectors are evaluated at each integration point $A_{T,q}$, $T\in\mathcal T_h$, $q\in\{1,2,\ldots, n_q\}$. So, we shall consider $6\times n_{int}$ arrays \verb+S+ and \verb+E+ for the stress and strain components, respectively. In particular, the array \verb+S+ will be used in for storage the values of $T_k$. Similarly,  a $36\times n_{int}$ array \verb+DS+ will store for the corresponding values of $T_k^o$ at all integration points.

The assembly of \verb+K+, \verb+F+, \verb+E+, \verb+S+ and \verb+DS+ is a subject of the next two sections. Section \ref{sec_assembly} deals with the elastic stiffness matrix, while Section \ref{sec_assembly_const} is focused on the assembly of the tangent stiffness matrix. From now on, we shall write \verb+K_elast+ and \verb+K_tangent+ to emphasize a type of the stiffness matrix.


\section{Assembly of elastic stiffness matrix}\label{sec_assembly}

The construction of the elastic stiffness matrix \verb+K_elast+ from  \eqref{K_elast-split} is realized by the function
\begin{verbatim}
 K_elast=elastic_stiffness_matrix(ELEM, COORD, shear, bulk, DHatP1, DHatP2, DHatP3, WF)
\end{verbatim}
within our codes. This function has the following input data:
\verb+ELEM+ is an $n_p\times n_e$ array that contains indices $1,2,\ldots, n_n$ of nodes belonging to each element and \verb+COORD+ is a $3\times n_n$ array containing coordinates of the nodes. We construct these standard arrays for specific geometries by vectorized procedures. Nevertheless, one can import them from a mesh generator. Further, \verb+shear+ and \verb+bulk+ denote $1\times n_{int}$ arrays representing values of the shear ($G$) and bulk ($K$) moduli at each integration points, respectively. These values are usually computed from the Young modulus and the Poisson ratio, see Section \ref{sec_example}. The $n_p\times n_q$ arrays \verb+DHatP1+, \verb+DHatP2+, \verb+DHatP3+ store the basis functions gradient values 
$$\frac{\partial\hat\Phi_p(\hat A_q)}{\partial\xi_1}, \frac{\partial\hat\Phi_p(\hat A_q)}{\partial\xi_2}, \frac{\partial\hat\Phi_p(\hat A_q)}{\partial\xi_3}, \qquad p=1,\ldots,n_p, \quad q=1,\ldots,n_q, $$ respectively. \verb+WF+ denotes a $1\times n_q$ array of the weight coefficients defining a numerical quadrature. The arrays \verb+WF+, \verb+DHatP1+, \verb+DHatP2+ and \verb+DHatP3+ are obtained by the following functions:
\begin{verbatim}
 [Xi, WF] = quadrature_volume(elem_type);
 [HatP,DHatP1,DHatP2,DHatP3] = local_basis_volume(elem_type, Xi);
\end{verbatim}         
where \verb+elem_type+ specifies P1, P2, Q1 or Q2 finite elements. Additionally, arrays \verb+Xi+, \verb+HatP+ contain local coordinates of $\hat A_q$ and the basis functions values 
$$\hat\Phi_p(\hat A_q), \qquad p=1,\ldots,n_p, \quad q=1,\ldots,n_q,$$ respectively. A quadrature rule is predefined for any element, see the Appendix. Nevertheless, one can easily change the rules within the function \verb+quadrature_volume+.

\subsection{Jacobian, its determinant and inverse, derivatives of local basis functions}
The Jacobian $J_T(\mbf\xi)$ from \eqref{Jacobian} needs to be evaluated at each integration point. Its components are stored in $1\times n_{int}$ arrays denoted as \verb+J11+, \verb+J12+,$\ldots$, \verb+J33+. These components are computed using the arrays \verb+DHatP1+, \verb+DHatP2+, \verb+DHatP3+, \verb+COORD+ and \verb+ELEM+. By a suitable replication of \verb+DHatP1+, \verb+DHatP2+, \verb+DHatP3+, we obtain the following $n_p\times n_{int}$ arrays:
\begin{verbatim}
 DHatPhi1=repmat(DHatP1,1,n_e); DHatPhi2=repmat(DHatP2,1,n_e); DHatPhi3=repmat(DHatP3,1,n_e);
\end{verbatim}
Further, from \verb+COORD+ and \verb+ELEM+, we derive $n_p\times n_q$ arrays \verb+COORDint1+, \verb+COORDint2+, \verb+COORDint3+ containing the first, second, and third coordinates (x, y and z components) of $n_p$ nodes that define an element containing a particular integration point:
\begin{verbatim}
 COORDe1=reshape(COORD(1,ELEM(:)),n_p,n_e); COORDint1=kron(COORDe1,ones(1,n_q));
 COORDe2=reshape(COORD(2,ELEM(:)),n_p,n_e); COORDint2=kron(COORDe2,ones(1,n_q));
 COORDe3=reshape(COORD(3,ELEM(:)),n_p,n_e); COORDint3=kron(COORDe3,ones(1,n_q));
\end{verbatim}
Now, one can easily compute the components of Jacobians at integration points:
\begin{verbatim}
 J11=sum(COORDint1.*DHatPhi1); J12=sum(COORDint2.*DHatPhi1); J13=sum(COORDint3.*DHatPhi1);
 J21=sum(COORDint1.*DHatPhi2); J22=sum(COORDint2.*DHatPhi2); J23=sum(COORDint3.*DHatPhi2);
 J31=sum(COORDint1.*DHatPhi3); J32=sum(COORDint2.*DHatPhi3); J33=sum(COORDint3.*DHatPhi3);
\end{verbatim}

Let \verb+DET+ and \verb+Jinv11+, $\ldots$, \verb+Jinv33+ be $1\times n_{int}$ arrays representing the determinant and the components of the inverse matrix to Jacobian. These arrays can be found by the following commands:
\begin{verbatim}
 DET = J11.*(J22.*J33-J23.*J32) - J12.*(J21.*J33-J23.*J31) + J13.*(J21.*J32-J22.*J31);
 Jinv11 =  (J22.*J33-J23.*J32)./DET; Jinv12 = -(J12.*J33-J13.*J32)./DET; Jinv13 =  (J12.*J23-J13.*J22)./DET; 
 Jinv21 = -(J21.*J33-J23.*J31)./DET; Jinv22 =  (J11.*J33-J13.*J31)./DET; Jinv23 = -(J11.*J23-J13.*J21)./DET; 
 Jinv31 =  (J21.*J32-J22.*J31)./DET; Jinv32 = -(J11.*J32-J12.*J31)./DET; Jinv33 =  (J11.*J22-J12.*J21)./DET; 
\end{verbatim}
According to (\ref{basis_derivative}), we evaluate 
$$\frac{\partial \Phi_{T,p}(A_{T,q})}{\partial x_1}, \frac{\partial \Phi_{T,p}(A_{T,q})}{\partial x_2}, \frac{\partial \Phi_{T,p}(A_{T,q})}{\partial x_3} \qquad T\in\mathcal T_h, \quad p=1,\ldots,n_p, \quad q=1,\ldots,n_q,$$ 
and store these values into $n_p\times n_{int}$ arrays \verb+DPhi1+, \verb+DPhi2+, \verb+DPhi3+:
\begin{verbatim}
 DPhi1 = repmat(Jinv11,n_p,1).*DHatPhi1 + repmat(Jinv12,n_p,1).*DHatPhi2 + repmat(Jinv13,n_p,1).*DHatPhi3;
 DPhi2 = repmat(Jinv21,n_p,1).*DHatPhi1 + repmat(Jinv22,n_p,1).*DHatPhi2 + repmat(Jinv23,n_p,1).*DHatPhi3;
 DPhi3 = repmat(Jinv31,n_p,1).*DHatPhi1 + repmat(Jinv32,n_p,1).*DHatPhi2 + repmat(Jinv33,n_p,1).*DHatPhi3;
\end{verbatim}

\subsection{Strain-displacement relation}

To represent a relation between the strain array \verb+E+ and the displacement array \verb+U+, we shall construct a $6n_{int}\times 3n_n$ array \verb+B+. By using this array, the strain-displacement relation can be written by the following command:
\begin{verbatim}
 E = reshape(B*U(:),6,n_int); 
\end{verbatim}
The array \verb+B+ is a large and sparse matrix, therefore its construction will be done by using the command \verb+sparse+,
\begin{verbatim}
 B = sparse(iB(:),jB(:),vB(:), 6*n_int,3*n_n);
\end{verbatim}
Here, \verb+vB+, \verb+iB+, \verb+jB+  are $18n_p\times n_{int}$ arrays containing non-zero values of \verb+B+ and the corresponding $i$-th and $j$-th indices, respectively. The size $18n_p=6*3n_p$ follows from the local strain-displacement relation defined at each quadrature point. It is well-known that the strain-displacement relation at a point $A_{T,q}$, $T\in\mathcal T_h$ and $q\in\{1,2,\ldots, r\}$ can be written as follows:
\begin{equation}
\left(
\begin{array}{c}
\varepsilon_{11}(A_{T,q})\\
\varepsilon_{22}(A_{T,q})\\
\varepsilon_{33}(A_{T,q})\\
2\varepsilon_{12}(A_{T,q})\\
2\varepsilon_{23}(A_{T,q})\\
2\varepsilon_{13}(A_{T,q})
\end{array}
\right)=
\left(
\begin{array}{ccccccc}
D_1\Phi^q_{T,1} & 0 & 0 & \ldots & D_1\Phi^q_{T,n_p} & 0 & 0 \\
0 & D_2\Phi^q_{T,1} & 0 & \ldots & 0 & D_2\Phi^q_{T,n_p} & 0 \\
0 & 0 & D_3\Phi^q_{T,1} & \ldots & 0 & 0 & D_3\Phi^q_{T,n_p} \\
D_2\Phi^q_{T,1} & D_1\Phi^q_{T,1} & 0 & \ldots & D_2\Phi^q_{T,n_p} & D_1\Phi^q_{T,n_p} & 0 \\
0 & D_3\Phi^q_{T,1} & D_2\Phi^q_{T,1} & \ldots & 0 & D_3\Phi^q_{T,n_p} & D_2\Phi^q_{T,n_p} \\
D_3\Phi^q_{T,1} & 0 & D_1\Phi^q_{T,1} & \ldots & D_3\Phi^q_{T,n_p} & 0 & D_1\Phi^q_{T,n_p}
\end{array}
\right)
\left(
\begin{array}{c}
U_1(N_{T,1})\\
U_2(N_{T,1})\\
U_3(N_{T,1})\\
\vdots\\
U_1(N_{T,n_p})\\
U_2(N_{T,n_p})\\
U_3(N_{T,n_p})
\end{array}
\right),
\label{local_strain-displ}
\end{equation}
where $D_i\Phi^q_{T,p}=\partial\Phi_{T,p}(A_{t,q})/\partial x_i$, $p=1,2,\ldots,n_p$, $i=1,2,3$. Each column of the array \verb+vB+ contains components of the matrix from (\ref{local_strain-displ}). In particular, we arrive at
\begin{verbatim}
 n_b=18*n_p; vB=zeros(n_b,n_int);
 vB(1:18:n_b-17,:)=DPhi1; vB(10:18:n_b- 8,:)=DPhi1; vB(18:18:n_b  ,:)=DPhi1;
 vB(4:18:n_b-14,:)=DPhi2; vB( 8:18:n_b-10,:)=DPhi2; vB(17:18:n_b-1,:)=DPhi2;
 vB(6:18:n_b-12,:)=DPhi3; vB(11:18:n_b- 7,:)=DPhi3; vB(15:18:n_b-3,:)=DPhi3;
\end{verbatim}
The arrays \verb+iB+ and \verb+jB+ of indices can be derived as follows:
\begin{verbatim}
 AUX=reshape(1:6*n_int, 6,n_int); 
 iB=repmat(AUX, 3*n_p,1);
 
 AUX1=[1;1;1]*(1:n_p); AUX2 = [2;1;0]*ones(1,n_p); AUX3=3*ELEM((AUX1(:))',:)-kron(ones(1,n_e),AUX2(:)); 
 jB=kron(AUX3,ones(6,n_q));
\end{verbatim}

\subsection{Elastic stiffness matrix}
\label{subsec_elast_assembly}

In elasticity, we have $T_k^o=\mathbb C=2G\mathbb I_D+K\mathbb I_V$ as follows from (\ref{elastic_operator}). To store these values at all integration points, we use the $36\times n_{int}$ array \verb+DS+ mentioned above:
\begin{verbatim}
 IOTA=[1;1;1;0;0;0]; VOL=IOTA*IOTA'; DEV=diag([1,1,1,1/2,1/2,1/2])-VOL/3; 
 DS=2*DEV(:)*shear+VOL(:)*bulk;      
\end{verbatim}
Here, \verb+IOTA+, \verb+VOL+, \verb+DEV+ denote the MATLAB counterparts of the tensors $\mbf I$, $\mathbb I_V=\mbf I\otimes\mbf I$, and $\mathbb I_D$, respectively. It is important to note that the definition of \verb+DEV+ enables to transform a strain-type tensor to a stress-type tensor. 

Further, we shall need the following $1\times n_{int}$ array
\begin{verbatim}
 WEIGHT = abs(DET).*repmat(WF, 1,n_e); 
\end{verbatim}
From the arrays \verb+DS+ and \verb+WEIGHT+, we arrive at the following $6n_{int}\times 6n_{int}$ block diagonal (sparse) matrix \verb+D_elast+:
\begin{verbatim}
 AUX=reshape(1:6*n_int,6,n_int);  
 iD=repmat(AUX,6,1); jD=kron(AUX,ones(6,1)); vD=DS.*repmat(WEIGHT,36,1);
 D_elast=sparse(iD,jD,vD,6*n_int,6*n_int );
\end{verbatim}
The elastic stiffness matrix can be assembled similarly as in formula (\ref{formula2}). 
The assembly of the array \verb+K_elast+ reads as
\begin{verbatim}
 K_elast=B'*D_elast*B; 
\end{verbatim}

\section{Assembly of the tangent stiffness matrix and vector of internal forces}
\label{sec_assembly_const}

The tangent stiffness matrix \verb+K_tangent+ based on \eqref{K-split} and the vector \verb+F+ of internal forces are updated in each time step and in each Newton iteration within the loading process. To this end, it suffices to update only the arrays \verb+S+, \verb+DS+ and use the arrays \verb+B+, \verb+iD+, \verb+jD+, \verb+WEIGHT+, \verb+K_elast+ and \verb+D_elast+ which are the output data from the function \verb+elastic_stiffness_matrix+. This is the main advantage of the presented assembly. 
The assembly of \verb+K_tangent+ and \verb+F+ read as
\begin{verbatim}
 vD = repmat(WEIGHT,9,1).*DS ;
 D_tangent = sparse( iD(:),jD(:),vD(:), 6*n_int,6*n_int ) ;   
 K_tangent = K_elast+B'*(D_tangent-D_elast)*B;   
 F = B'*reshape(S.*repmat(WEIGHT,6,1), 6*n_int,1);
\end{verbatim}
and the arrays \verb+S+ and \verb+DS+ are created by the function
\begin{verbatim}
 constitutive_problem
\end{verbatim}
described in Sections \ref{subsec_VM_assembly} and \ref{subsec_DP_assembly} von von Mises and Drucker-Prager yield criteria. 

\subsection{Von Mises yield criterion and kinematic hardening}
\label{subsec_VM_assembly}

According to Section \ref{subsec_kinematic}, the input data to the function \verb+constitutive_problem+ are:
\begin{verbatim}
 E, Ep_prev, Hard_prev, shear, bulk, a , Y.
\end{verbatim}
Here, \verb+E = reshape(B*U(:),6,n_int)+, \verb+Ep_prev+ and \verb+Hard_prev+ are $6\times n_{int}$ arrays representing the strain tensor $\mbf\varepsilon_k$ at current time step $k$, the plastic strain $\mbf\varepsilon_{k-1}^p$ and the kinematic hardening $\mbf\beta_{k-1}$ from the previous time step, respectively. The remaining $1\times n_{int}$ input arrays store the material parameters $G$, $K$, $a$ and $Y$ at all integration points. 

First, we compute the MATLAB counterparts to $\mbf \varepsilon_k-\mbf \varepsilon^p_{k-1}$, $\mbf \sigma_k^{tr}$, $\mbf s_k^{tr}$, and $|\mbf s_k^{tr}|$:
\begin{verbatim}
 E_tr=E-Ep_prev;                                                                   % size(E_tr)=(6,n_int)
 S_tr=2*repmat(shear,6,1).*(DEV*E_tr)+repmat(bulk,6,1).*(VOL*E_tr);                % size(S_tr)=(6,n_int)      
 SD_tr=DEV*(2*repmat(shear,6,1).*E_tr)-Hard_prev;                                  % size(SD_tr)=(6,n_int)
 norm_SD=sqrt(sum(SD_tr(1:3,:).*SD_tr(1:3,:))+2*sum(SD_tr(4:6,:).*SD_tr(4:6,:)));  % size(norm_SD)=(1,n_int)
\end{verbatim}
In order to distinguish integration points with elastic and plastic behaviour, we define a $1\times n_{int}$ array \verb+CRIT+ representing the yield criterion and the corresponding logical $1\times n_{int}$ array \verb+IND_p+ which indicates integration points with plastic behaviour:
\begin{verbatim}
 CRIT=norm_SD-Y;  IND_p=CRIT>0;  
\end{verbatim}
The elastic prediction yields
\begin{verbatim}
 S=S_tr; DS=2*DEV(:)*shear+VOL(:)*bulk;
\end{verbatim}
We apply the plastic correction at the integration points with the plastic response according to the formulas (\ref{T_kinematic}) and (\ref{DT_kinematic}):
\begin{verbatim}
 Nhat=SD_tr(:,IND_p)./repmat(norm_SD(IND_p),6,1);
 denom = 2*shear(IND_p)+a(IND_p);  lambda=CRIT(IND_p)./denom; 
 S(:,IND_p)=S(:,IND_p)-repmat(2*shear(IND_p).*lambda,6,1).*N_hat;
\end{verbatim}
Here, the arrays \verb+Nhat+, \verb+lambda+ represent the function $\mbf n^{tr}_k$ and the plastic multipliers. In order to update \verb+DS+, we introduce two auxiliary $36\times n_{int}$ arrays \verb+ID+ and \verb+NNhat+ representing the terms $\mathbb I_D$ and $\mbf n_k^{tr}\otimes\mbf n_k^{tr}$ in (\ref{DT_kinematic}), respectively: 
\begin{verbatim}
 ID=DEV(:)*ones(1,length(lambda));
 NNhat=repmat(N_hat,6,1).*kron(N_hat,ones(6,1));
 const=((2*shear(IND_p)).^2)./denom;
 DS(:,IND_p)=DS(:,IND_p)-repmat(const,36,1).*ID+...
                repmat((const.*Y(IND_p))./norm_SD(IND_p),36,1).*(ID-NN_hat);
\end{verbatim}
Let us complete that the function \verb+constitutive_problem+ also contains other output data including the updated plastic strain and the hardening according to formula (\ref{T_other}).

\subsection{Drucker-Prager yield criterion and perfect plasticity}
\label{subsec_DP_assembly}

The constitutive assembly for the Drucker-Prager yield criterion and perfect plasticity is based on formulas from Section \ref{subsec_DP}. The input data to the function \verb+constitutive_problem+ are now \verb+E+, \verb+E_prev+, \verb+shear+, \verb+bulk+ and $1\times n_{int}$ arrays \verb+eta+ and \verb+c+ representing plastic material parameters. The corresponding MATLAB code is summarized below. Since its structure is similar to Section \ref{subsec_VM_assembly}, we skip some comments to this code. We only emphasize that two decision criteria are used unlike the von Mises model. To this end, we introduce $1\times n_{int}$ logical arrays \verb+IND_s+ and \verb+IND_a+ indicating integration points where the return to the smooth portion and to the apex of the yield surface happen, respectively.

\begin{verbatim}
  E_tr=E-Ep_prev;                          
  S_tr=2*repmat(shear,6,1).*(DEV*E_tr)+repmat(bulk,6,1).*(VOL*E_tr);   
  dev_E=DEV*E_tr;                       % deviatoric part of E_tr
  norm_E=sqrt(max(0,sum(E_tr.*dev_E))); % norm of the deviatoric strain
  rho_tr=2*shear.*norm_E;               % \varrho^{tr}
  p_tr=bulk.*(IOTA'*E_tr);              % trial volumetric stress
  
  denom_a= bulk.*(eta.^2);
  denom_s=shear+denom_a;
  CRIT1= rho_tr/sqrt(2) + eta.*p_tr - c ; 
  CRIT2= eta.*p_tr - denom_a.*rho_tr./(shear*sqrt(2)) - c ;
  
  IND_s = (CRIT1>0)&(CRIT2<=0);   % logical array for the return to the smooth portion 
  IND_a = (CRIT1>0)&(CRIT2>0);    % logical array for the return to the apex 

  S=S_tr; DS=2*DEV(:)*shear+VOL(:)*bulk; 

  lambda_s=CRIT1(IND_s)./denom_s(IND_s); 
  n_smooth=length(lambda_s);  
  lambda_a=(eta(IND_a).*p_tr(IND_a)-c(IND_a))./denom_a(IND_a);
  n_apex=length(lambda_a);
  
  N_hat=dev_E(:,IND_s)./repmat(norm_E(IND_s),6,1);
  M_hat=repmat(sqrt(2)*shear(IND_s),6,1).*N_hat+IOTA*(bulk(IND_s).*eta(IND_s));
  S(:,IND_s)=S(:,IND_s)-repmat(lambda_s,6,1).*M_hat;
  S(:,IND_a)=IOTA*(c(IND_a)./eta(IND_a));

  ID=DEV(:)*ones(1,n_smooth);
  NN_hat=repmat(N_hat,6,1).*kron(N_hat,ones(6,1));
  MM_hat=repmat(M_hat,6,1).*kron(M_hat,ones(6,1));
  DS(:,IND_s)=DS(:,IND_s)-...
      repmat(2*sqrt(2)*(shear(IND_s).^2).*lambda_s./rho_tr(IND_s),36,1).*(ID-NN_hat)-...
      MM_hat./repmat(denom_s(IND_s),36,1);
  DS(:,IND_a)=zeros(36,n_apex);
\end{verbatim}
Let us complete that the function \verb+constitutive_problem+ also contains other output  updating the plastic strain according to formulas from Section \ref{subsec_DP}.


\section{Computational Examples}
\label{sec_example}
Particular examples are introduced in Sections \ref{subsec_experiment_el}-\ref{subsec_experiment_DP}. The performance of the codes was tested with MATLAB 8.5.0.197613 (R2015a) on a computer with 64 Intel(R) Xeon(R) CPU E5-2640 v3 processors running at 2.60GHz, number of processors 2, 256 GB RAM, and 7.5 TB harddisk memory. Calculation times are given in seconds. Our codes  are  available for download and testing at \cite{software}. They contain implementation of selected 2D (plain strain) and 3D problems from elasticity and elastoplasticity on predefined domains. One can choose P1, P2, Q1 or Q2 elements (see Appendix) and different levels of mesh density. Regular meshes are considered and constructed using fully vectorized procedures.

\subsection{Assembly of elastic stiffness matrices and comparison with other codes}
\label{subsec_experiment_el}
\begin{figure}
  \begin{minipage}[t]{0.47\textwidth}
        \centering
        \begin{picture}(140,140)

            {\thicklines
            \put(70,20){\line(1,0){50}}
            \put(20,120){\line(1,0){100}}
            \put(20,70){\line(0,1){50}}
            \put(120,20){\line(0,1){100}}
            \put(70,20){\line(0,1){50}}
            \put(20,70){\line(1,0){50}}
            }
            \multiput(15,72)(0,5){10}{\circle{4}}
            \multiput(72,15)(5,0){10}{\circle{4}}
            \multiput(22,122)(8,0){13}{\vector(0,1){10}}
            \multiput(72,82)(8,0){3}{\vector(0,1){7}}
            \put(75,25){\vector(1,0){40}}
            \put(20,20){\vector(1,0){20}}
            \put(20,20){\vector(0,1){20}}

            \put(70,134){\makebox(0,0)[b]{$f_T$}}
            \put(80,92){\makebox(0,0)[b]{$f_V$}}
            \put(95,27){\makebox(0,0)[b]{$u_{D}$}}
            \put(40,23){\makebox(0,0)[b]{$x_1$}}
            \put(23,40){\makebox(0,0)[l]{$x_2$}}
            \put(45,0){\makebox(0,0)[t]{5}}
            \put(95,0){\makebox(0,0)[t]{5}}
            \put(132,70){\makebox(0,0)[l]{10}}
            \put(20,4){\line(1,0){100}}
            \put(130,20){\line(0,1){100}}
            \put(20,2){\line(0,1){4}}
            \put(70,2){\line(0,1){4}}
            \put(120,2){\line(0,1){4}}
            \put(128,20){\line(1,0){4}}
            \put(128,120){\line(1,0){4}}

        \end{picture}
  \end{minipage}
  \hfill
  \begin{minipage}[t]{0.47\textwidth}
    \center
    \includegraphics[width=0.9\textwidth]{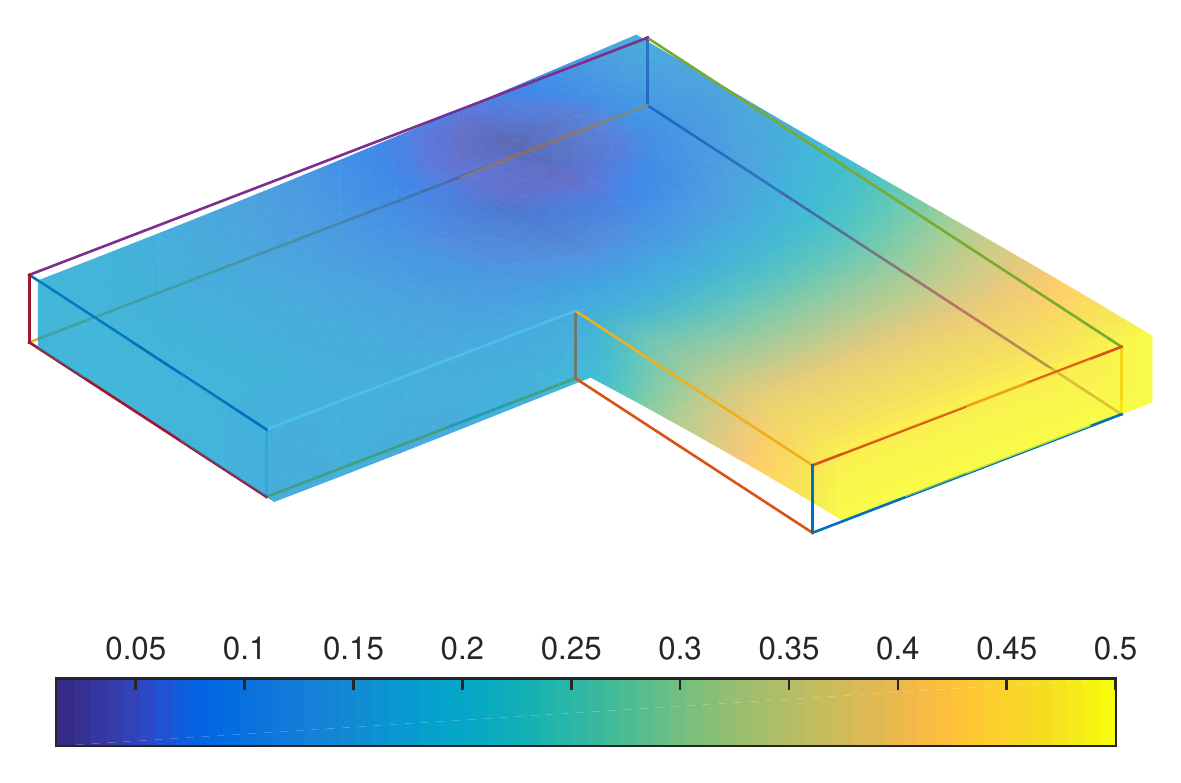}
  \end{minipage}
  \caption{Simplified 2D geometry of the elastic problem (left). The real 3D geometry appears by extrusion in the $x_3$ direction. The corresponding total displacement field is displayed in a deformed configuration (right).}
  \label{fig_scheme}
\end{figure}
We consider a body that occupies the domain depicted in Figure \ref{fig_scheme} in $x_1-x_2$ plane. The corresponding 3D geometry appears by extrusion in $x_3$ direction. The size of the body in this direction is equal to one if the 3D problem is considered. It is assumed that prescribed forces are independent of $x_3$ direction to be 2D and 3D results similar. On the left and bottom sides of the depicted domain, the symmetry boundary conditions are prescribed, i.e., $\mbf u\cdot\mbf n=0$ where $\mbf n$ is a normal vector to the boundary. On the bottom, we also prescribe nonhomogeneous Dirichlet boundary condition $u_{D}=0.5$ in the direction $x_1$. Further, the constant traction of density $f_t=200$ is acting on the upper side in the normal direction and the constant volume force $F_V=1$ is prescribed in $x_2$ direction.  The material parameters are set as follows: $E = 206 900\,$ (Young's modulus) and $\nu = 0.29$ (Poisson's  ratio). The corresponding values of the bulk and shear moduli are computed standardly from $E$ and $\nu$: $$K=\frac{E}{3(1-2\nu)}, \quad G=\frac{E}{2(1+\nu)}.$$



\begin{figure}
\center
\begin{minipage}{0.45\textwidth}
\center
\includegraphics[width=\textwidth]{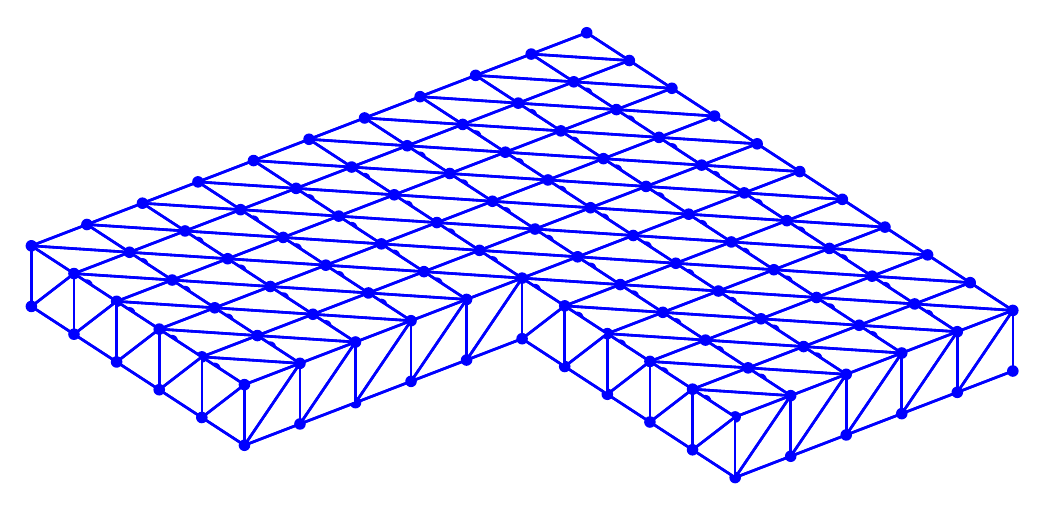} \\
$P1$ \\
\includegraphics[width=\textwidth]{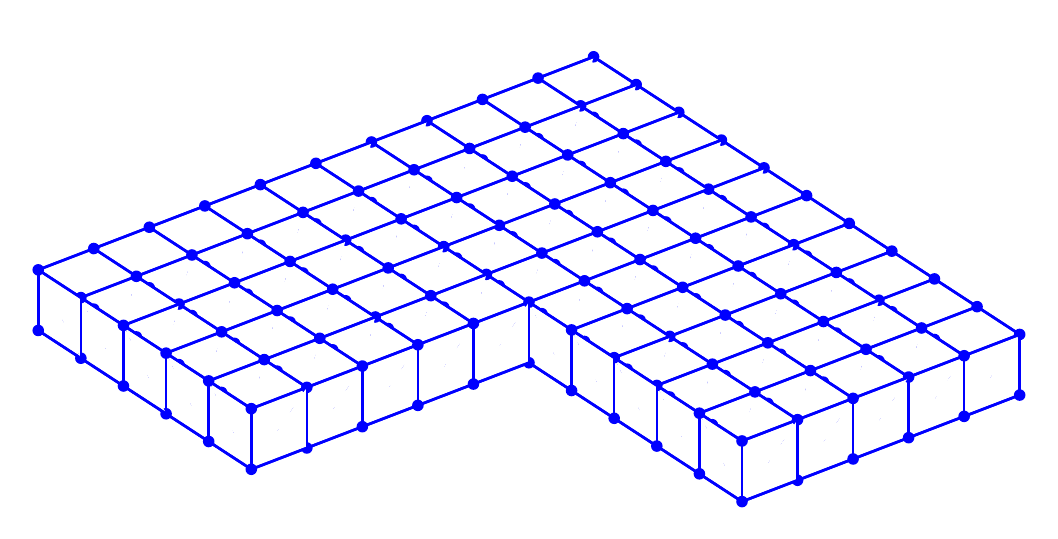} \\
$Q1$
\end{minipage}
\vspace{0.03\textwidth}
\begin{minipage}{0.45\textwidth}
\center
\includegraphics[width=\textwidth]{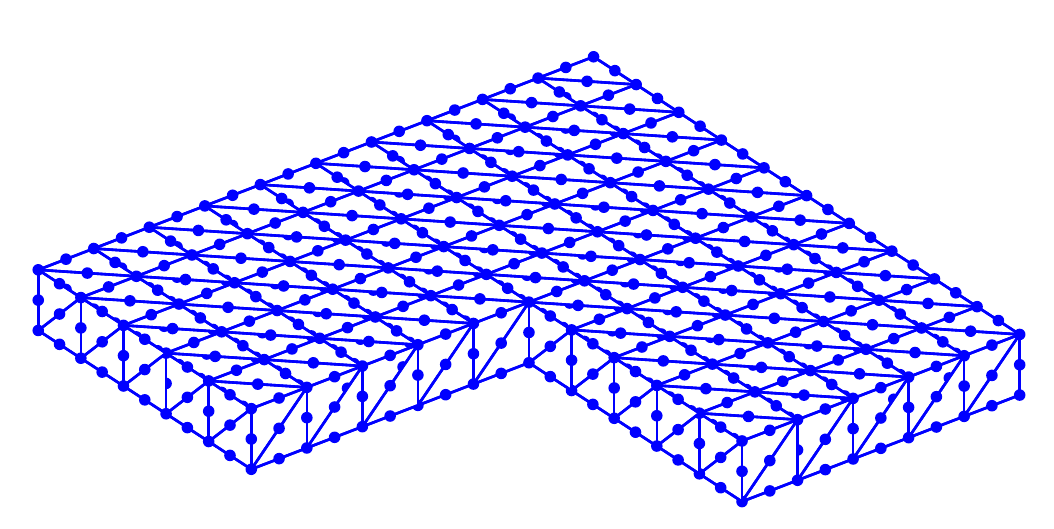} \\
$P2$ \\
\includegraphics[width=\textwidth]{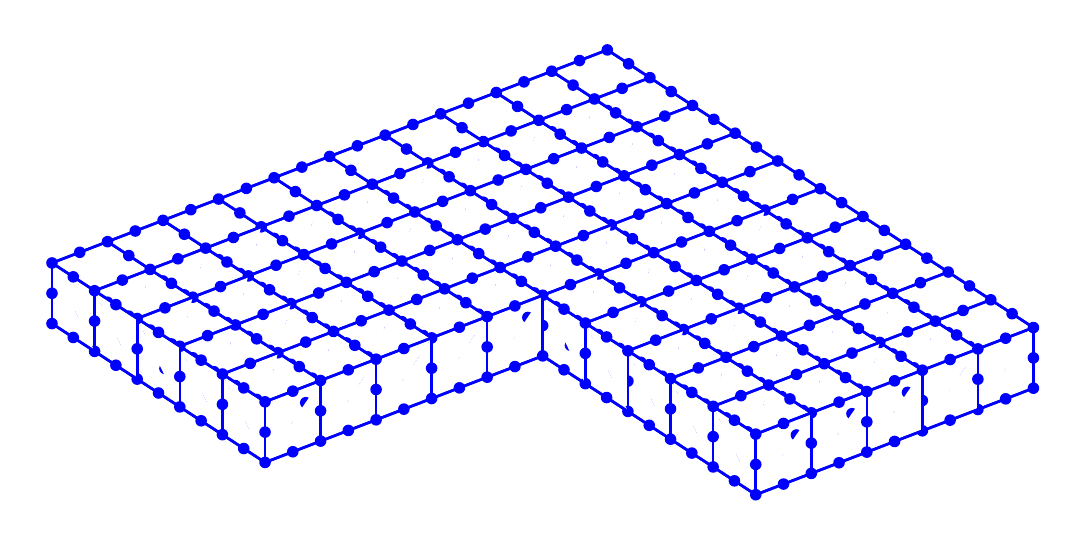} \\
$Q2$
\end{minipage}
\caption{A 3D domain discretized by P1, P2, Q1 and Q2 elements.} 
\label{fig:3D_elements}
\end{figure}

The corresponding codes are located in directory \verb+elasticity+ at \cite{software}. Stiffness matrices are generated for a sequence of uniformly refined meshes in 2D and 3D. Coarse (level 0) 3D meshes are displayed in Figure \ref{fig:3D_elements}. Results were obtained by the script
\begin{verbatim}
  elasticity_assembly_test;
\end{verbatim}
located in subdirectories \verb+elasticity_2D+ and \verb+elasticity_3D+. 
Corresponding assembly times are summarized in Tables \ref{ta:timesP1_2D_elasticity}, \ref{ta:timesP1_3D_elasticity}. 
\begin{table}[ht]
 \centering
\begin{tabular}{   r || r | r | r |  r | |  r  |  r | r  |  r   | r  |  r}
 &   \multicolumn{4}{c}{size of}  &  \multicolumn{4}{c}{assembly of}  \\
 level &  $\stiff^{P1}$ & $\stiff^{P2}$ & $\stiff^{Q1}$ & $\stiff^{Q2}$ & $\stiff^{P1}$ & $\stiff^{P2}$  & $\stiff^{Q1}$ & $\stiff^{Q2}$  \\
\hline
3 &         9922 &      39,042 &        9922 &      29,442 &   0.06 &      0.62 &    0.17 &    0.56 \\
4 &       39,042 &     154,882 &      39,042 &     116,482 &   0.21 &      2.40 &    0.58 &    1.89 \\
5 &      154,882 &     616,962 &     154,882 &     463,362 &   0.72 &      9.13 &    1.70 &    8.04 \\
6 &      616,962 &   2,462,722 &     616,962 &   1,848,322 &   2.75 &     34.58 &    6.70 &   32.19 \\
7 &    2,462,722 &   9,840,642 &   2,462,722 &   7,383,042 &  10.42 &    148.68 &   27.17 &  130.65 \\
8 &    9,840,642 &  39,342,082 &   9,840,642 &  29,511,682 &  39.87 &    608.70 &  112.40 &  545.08 \\
9 &   39,342,082 & 157,327,362 &  39,342,082 & 118,005,762 & 169.52 & 16,524.98 &  528.65 & 8284.21 \\
\end{tabular}
\caption{2D assembly of {\bf{elastic}} stiffness matrices for P1, P2, Q1, Q2 elements.}\label{ta:timesP1_2D_elasticity}
\vspace{0.5cm}
\begin{tabular}{   r || r | r | r |  r | |  r  |  r | r  |  r   | r  |  r}
 &   \multicolumn{4}{c}{size of}  &  \multicolumn{4}{c}{assembly of}  \\
 level &  $\stiff^{P1}$ & $\stiff^{P2}$ & $\stiff^{Q1}$ & $\stiff^{Q2}$ & $\stiff^{P1}$ & $\stiff^{P2}$  & $\stiff^{Q1}$ & $\stiff^{Q2}$  \\
\hline
1 &      3069 &    19,215 &      3069 &    10,875 &   0.14 &    1.55 &   0.19 &   1.51 \\
2 &    19,215 &   133,947 &    19,215 &    71,787 &   0.43 &   12.25 &   1.18 &  11.54 \\
3 &   133,947 &   995,571 &   133,947 &   516,531 &   3.70 &  102.07 &   9.78 &  95.87 \\
4 &   995,571 & 7,666,659 &   995,571 & 3,907,299 &  29.13 & 1040.36 &  78.38 & 849.60 \\
\end{tabular}
\caption{3D assembly of {\bf{elastic}} stiffness matrices for for P1, P2, Q1, Q2 elements.}
\label{ta:timesP1_3D_elasticity}
\end{table}
We observe (almost) optimal scalability: assembly times are linearly proportional to sizes of matrices. The solution of the elastic problem defined above is computed by the script
\begin{verbatim}
  elasticity_fem;
\end{verbatim}
located in subdirectories \verb+elasticity_2D+ and \verb+elasticity_3D+. 
Beside elastic stiffness matrix $\mbf K_{elast}$, the right-hand size vector $\mbf f$ is  assembled by fully vectorized procedures. The procedures are written for more general volume and surface forces that need not be only constant. The displacement $\mbf u$ is computed from the linear system of equations 
$$\mbf K_{elast} \mbf u = \mbf f$$ and displayed together with the deformed body, see Figure \ref{fig_scheme}. 
 

Performance comparison to the technique of Rahman and Valdman \cite{RaVa} for P1 elements is done by scripts
\begin{verbatim}
  comparison_assembly_P1_2D_elasticity;
  comparison_assembly_P1_3D_elasticity;
\end{verbatim}
located in the main directory 
and reported in Tables \ref{ta:timesP1_2D_elasticity_comparison}, \ref{ta:timesP1_3D_elasticity_comparison}. 
\begin{table}[ht]
\begin{minipage}{0.48\textwidth}\centering
\begin{tabular}{  r || r | r  |  r }
 &  size of  &  \multicolumn{2}{r}{assembly of}  \\
 level &  $\stiff^{P1}, \stiff^{P1}_{RV}$  &  $\stiff^{P1}$ & $\stiff^{P1}_{RV}$ \\
\hline
6  &      25,090 &   0.12 &    0.22 \\
7  &      99,330 &   0.53 &    1.01 \\
8  &     395,266 &   2.25 &    4.31 \\
9  &   1,576,962 &   9.54 &   16.23 \\
10 &   6,299,650 &  43.01 &   70.02 \\
11 &  25,182,210 & 180.15 &  308.68 \\
12 & 100,696,066 & 975.82 & 1536.63 \\
\end{tabular}
\caption{2D assembly of {\bf{elastic}} stiffness matrices for P1 elements. 
}
\label{ta:timesP1_2D_elasticity_comparison}
\end{minipage}
\hfill
\begin{minipage}{0.48\textwidth}\centering
\begin{tabular}{  r || r | r  |  r }
 &  size of  &  \multicolumn{2}{r}{assembly of}  \\
 level &  $\stiff^{P1}, \stiff^{P1}_{RV}$  &  $\stiff^{P1}$ & $\stiff^{P1}_{RV}$ \\
\hline
1 &       1029 &    0.03 &    0.05 \\
2 &       6591 &    0.19 &    0.42 \\
3 &     46,875 &    1.62 &    3.96 \\
4 &    352,947 &   11.36 &   30.55 \\
5 &  2,738,019 &  101.36 &  255.56 \\
6 & 21,567,171 & 1055.36 & 3742.85 \\
\end{tabular}
\caption{3D assembly of {\bf{elastic}} stiffness matrices for P1 elements. 
}
\label{ta:timesP1_3D_elasticity_comparison}
\end{minipage}
\end{table}
Our technique is about 2 times faster in 2D and 3 times faster in 3D. Another comparison with 3D techniques of \cite{ACFK02, Koko} can be run by the script 
 \begin{verbatim}
  comparison_fem_3D_elasticity;
\end{verbatim}
The original  assembly  of the function \verb+fem_lame3d+ of \cite{ACFK02} requires 6.98 seconds, the modification by our technique 0.25 seconds. This huge improvement is due to the fact that the assembly of \cite{ACFK02} is not vectorized. The original  assembly  of the function \verb+demo_elas+ of \cite{Koko} requires 0.33 seconds, the modification by our technique comparable 0.21 seconds.

\subsection{Assembly of plastic stiffness matrices for the von Mises yield criterion}
\label{subsec_experiment_VM}

\begin{figure}[htbp]
        \centering
        \begin{picture}(230,130)
          
            \put(0,60){\vector(1,0){230}}
            \put(10,0){\vector(0,1){130}}
            \put(10,60){\line(1,1){50}}
            \put(60,110){\line(1,-1){100}}
            \put(160,10){\line(1,1){50}}
    
            \put(13,130){\makebox(0,0)[l]{$f_T$}}
            \put(230,63){\makebox(0,0)[b]{$t$}}
            \put(60,56){\makebox(0,0)[t]{1}}
            \put(110,56){\makebox(0,0)[t]{2}}
            \put(160,56){\makebox(0,0)[t]{3}}
            \put(210,56){\makebox(0,0)[t]{4}}
            \put(6,10){\makebox(0,0)[r]{-200}}
            \put(6,110){\makebox(0,0)[r]{200}}
            \put(60,58){\line(0,1){4}}
            \put(160,58){\line(0,1){4}}
            \put(8,110){\line(1,0){4}}
            \put(8,10){\line(1,0){4}}

        \end{picture}
\caption{History of the traction force.}
\label{fig_traction}
\end{figure}
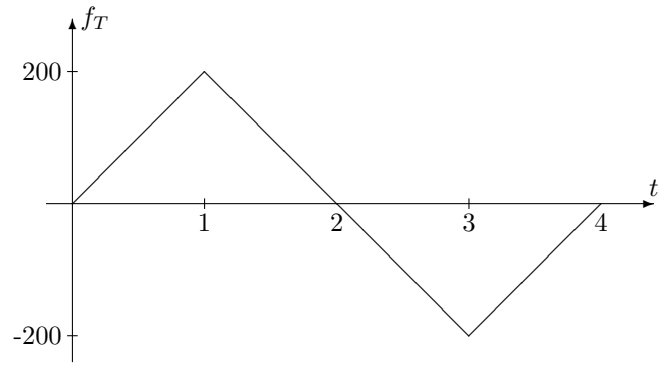

We consider the same geometry as in Section \ref{subsec_experiment_el}, see Figure \ref{fig_scheme}. Unlike elasticity, we do not consider the volume force and the nonhomogeneous Dirichlet boundary condition. The traction force prescribed on the upper side is now time dependent, see Figure \ref{fig_traction}. So the traction force is of the form $\zeta(t)f_{T,max}$, $t\in[0,4]$, where the scale $\zeta$ of external forces varies from $-1$ to $1$, and $f_{T, max}=200$. The inelastic material parameters are set as follows: $a=10,000$ and $Y=450\sqrt{2/3}$.

\begin{figure}
\center
\begin{minipage}{0.45\textwidth}
\center
\includegraphics[width=\textwidth]{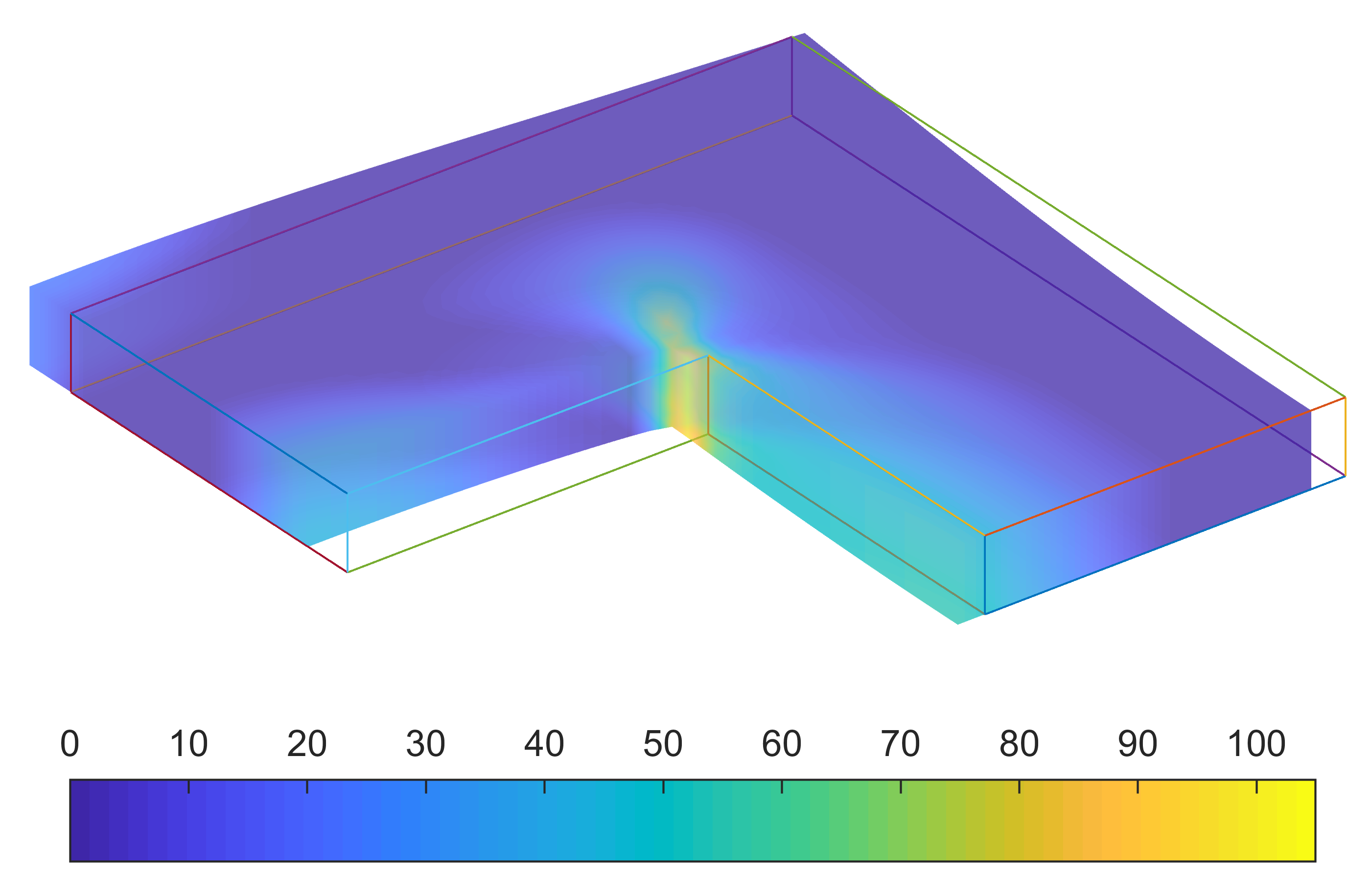} \\
$t=10$ \\
\hspace{0.2cm}
\includegraphics[width=\textwidth]{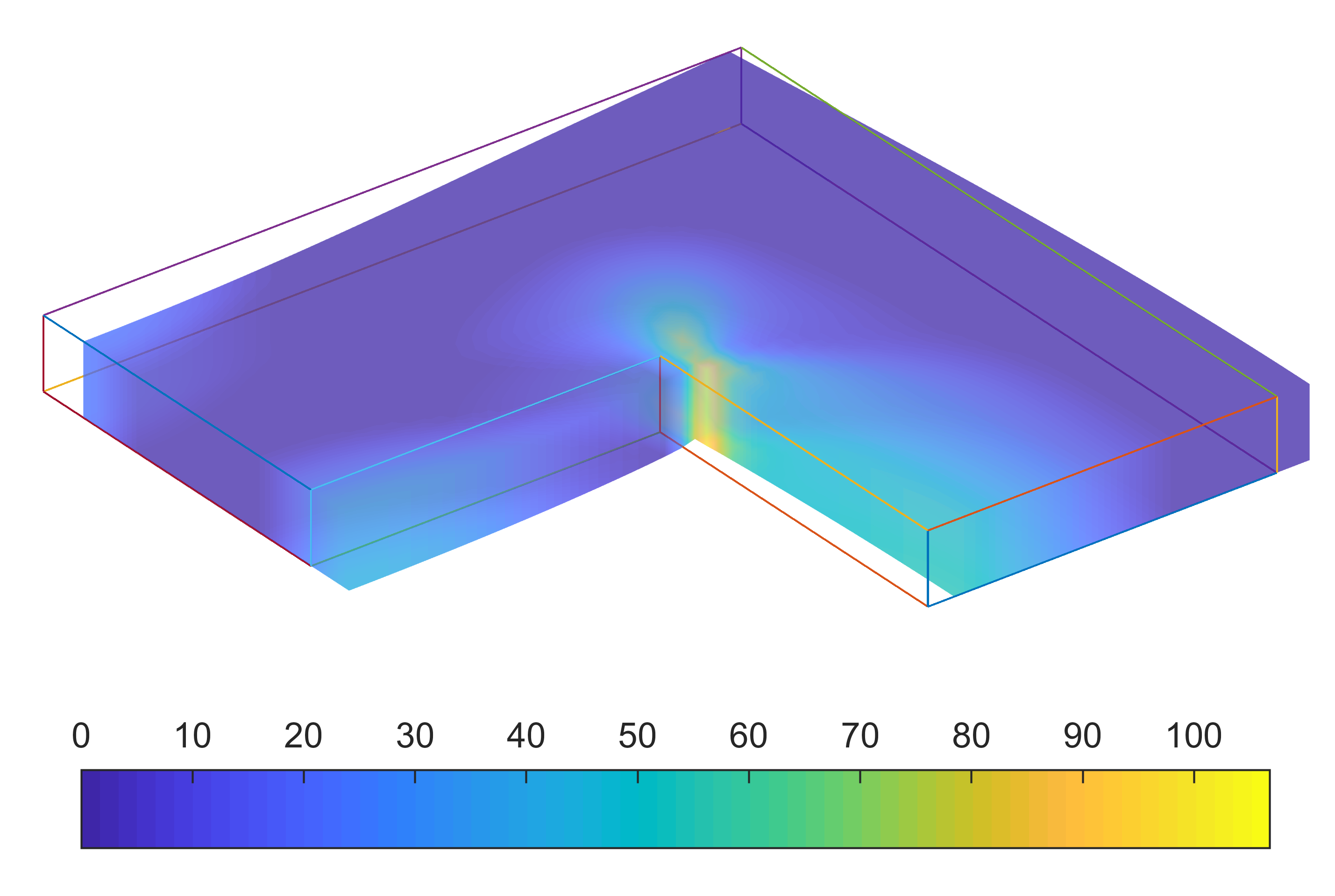} \\
$t=30$
\end{minipage}
\vspace{0.03\textwidth}
\begin{minipage}{0.45\textwidth}
\center
\includegraphics[width=\textwidth]{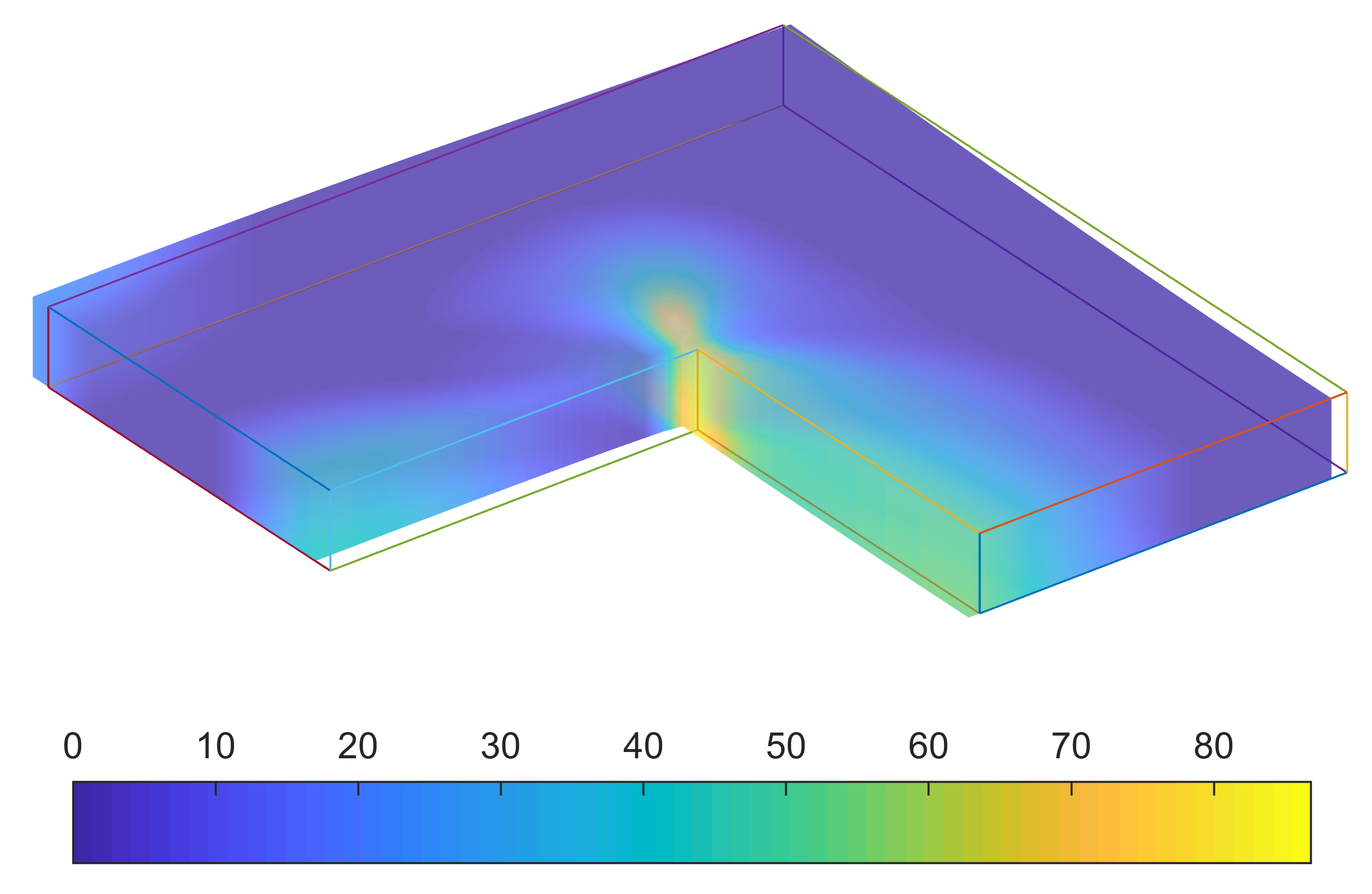} \\
$t=20$ \\
\hspace{0.2cm}
\includegraphics[width=\textwidth]{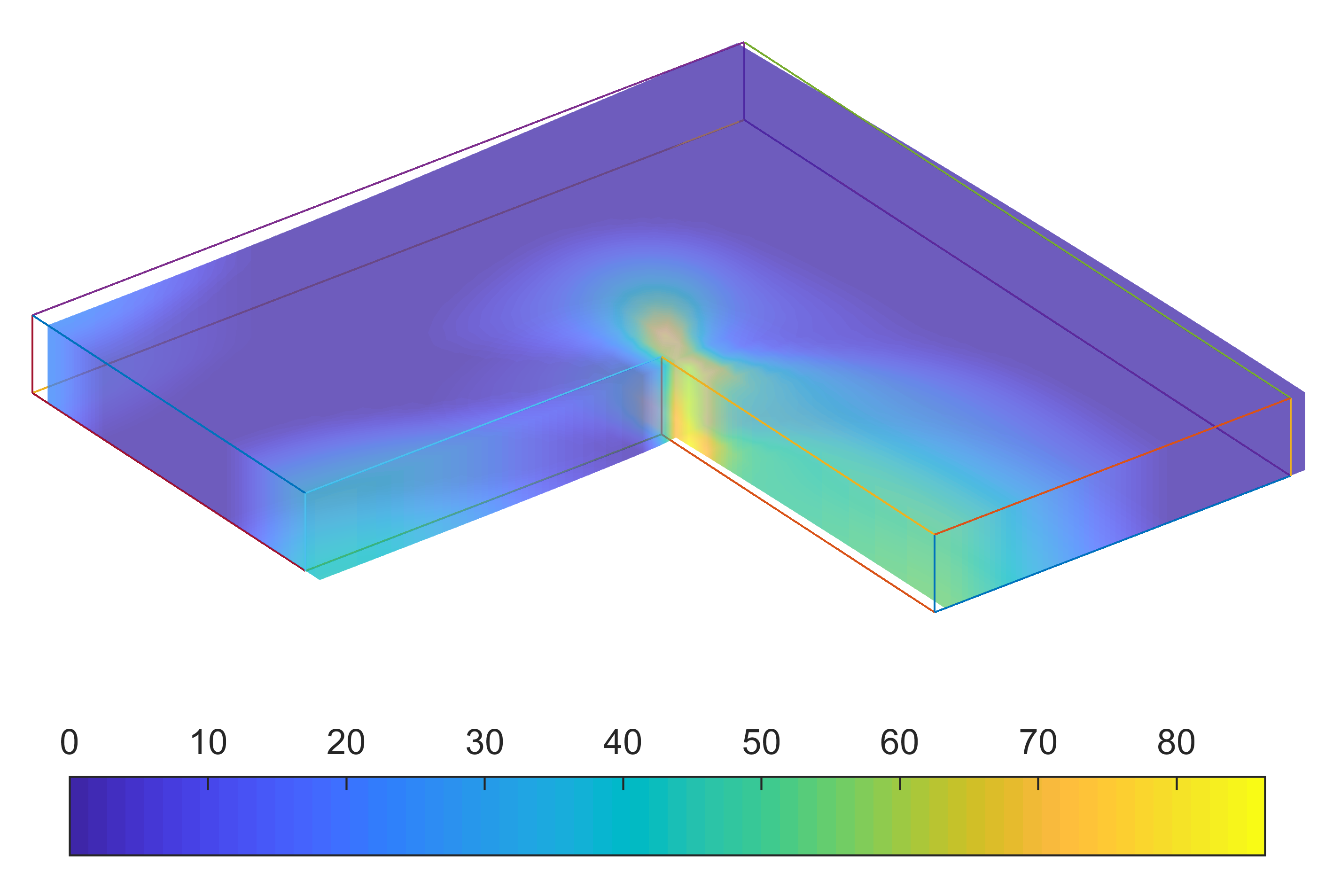} \\
$t=40$
\end{minipage}
\caption{Hardening fields at discrete times $10, 20, 30, 40$.}
\label{fig:hardening}
\vspace{1cm}
\includegraphics[width=0.45\textwidth]{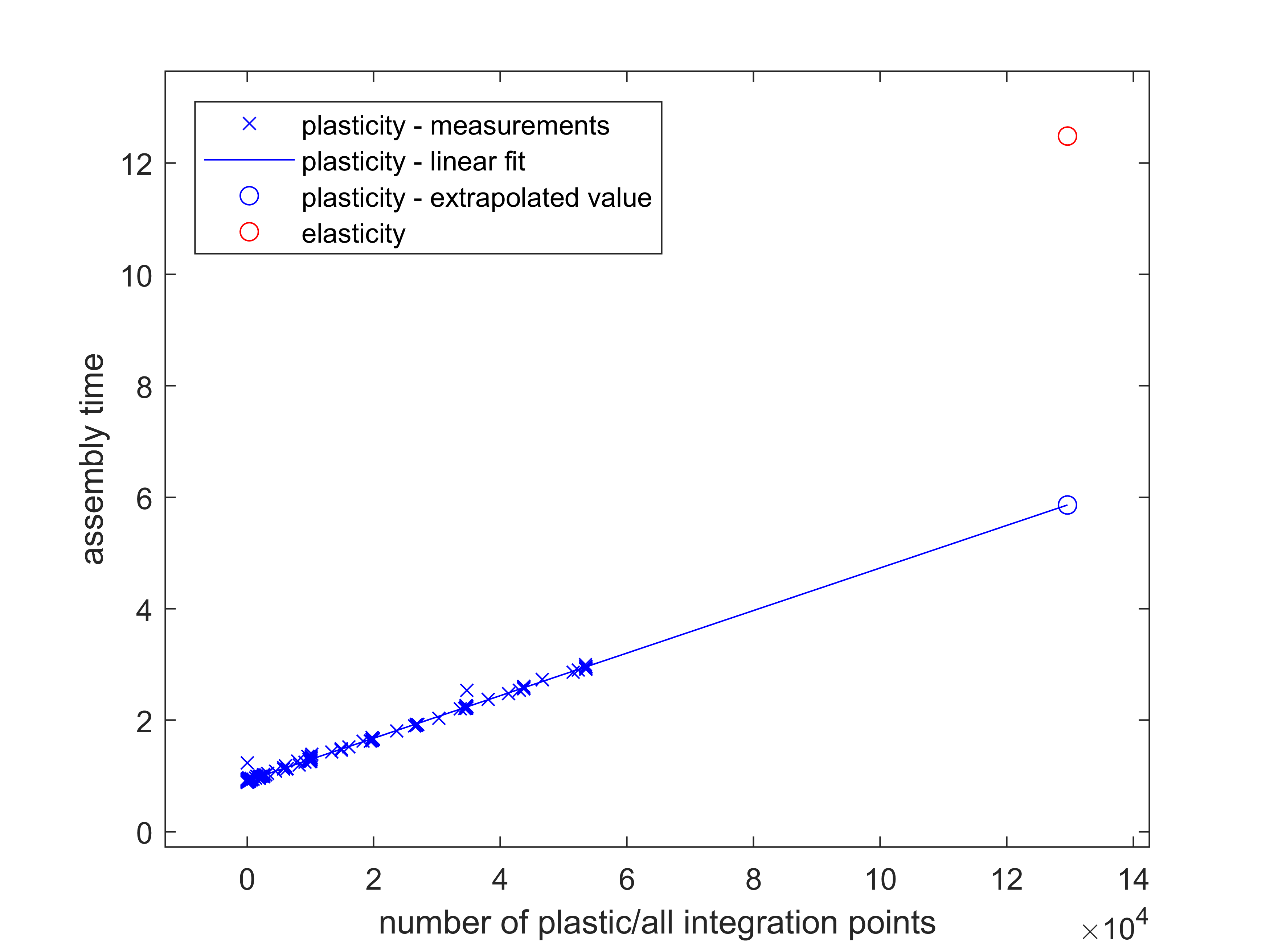}
\vspace{0.03\textwidth}
\includegraphics[width=0.45\textwidth]{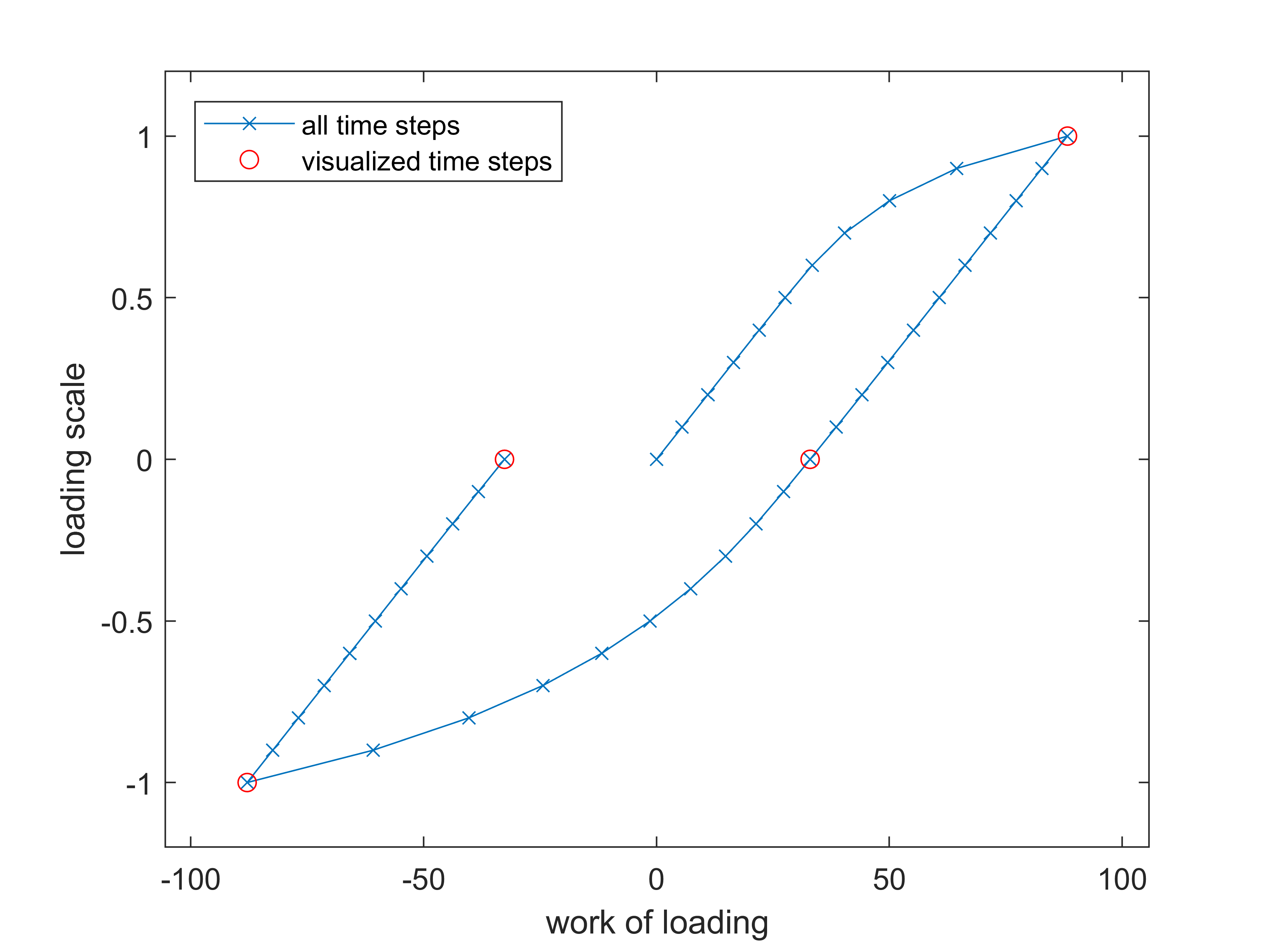}
\caption{Assembly times of tangential stiffness matrix versus number of plastic integration points (left) and a hysteresis curve (right).}
\label{fig:hysteresis}
\end{figure}

Elastoplasticity related codes are located in directory \verb+plasticity+ and a sequence of incremental steps using von Mises criterion is solved by the script 
\begin{verbatim}
  plasticity_VM_fem;
\end{verbatim}
located in subdirectories \verb+plasticity_VM_2D+ and \verb+plasticity_VM_3D+. In each step, few iterations of the semismooth Newton method are performed until the convergence is
reached. 

We visualize the results for Q2 elements and the mesh with 23,929 nodes (level 2). Figure \ref{fig:hardening} depicts hardening fields for Q2 elements  in selected time steps.
Assembly times of the tangential stiffness matrix in each Newton iteration are stored together with the number of integration points in the plastic regime. 
Results are displayed in Figure \ref{fig:hysteresis}. There is a linear relation between the assembly times and the numbers of integration points in the plastic regime. The elastic stiffness matrix is precomputed and its assembly time is not added to measured times. The extrapolated value suggests if all integration points are in the plastic regime, the assembly takes half the time of the elastic matrix assembly in the worst case. Additionally, we also show a hysteresis curve in the same figure as a relation of the scale $\zeta$ of external forces and the work of external forces. The work is computed as $\mbf f_{max}^T\mbf u_k$, where $\mbf f_{max}^T$ is a vector representing the maximal traction force $f_{T,max}$ and $\mbf u_k$ is a solution of discretized problem at $k$th time step.

\subsection{Assembly of plastic stiffness matrices for the Drucker-Prager yield criterion}
\label{subsec_experiment_DP}

An example of computations with Drucker-Prager criterion is available in the script 
\begin{verbatim}
  plasticity_DP_fem;
\end{verbatim}
located in subdirectories \verb+plasticity_DP_2D+ and \verb+plasticity_DP_3D+. It simulates a well know strip-footing benchmark \cite[Chapters 7,8]{NPO08} leading to bearing capacity (limit load) of a soil foundation. 

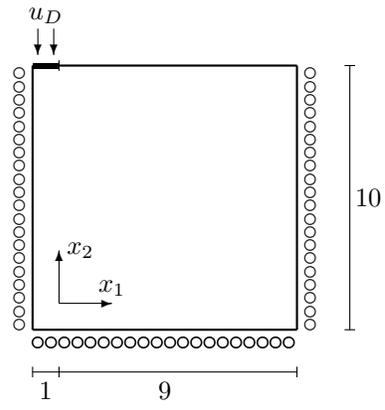
\begin{figure}[htbp]
        \centering
        \begin{picture}(150,140)

            {\thicklines
            \put(20,20){\line(1,0){100}}
            \put(20,120){\line(1,0){100}}
            \put(20,20){\line(0,1){100}}
            \put(120,20){\line(0,1){100}}
             }
            \multiput(15,22)(0,5){20}{\circle{4}}
            \multiput(22,15)(5,0){20}{\circle{4}}
            \multiput(125,22)(0,5){20}{\circle{4}}
            \multiput(22,15)(5,0){20}{\circle{4}}
            \multiput(22,134)(6,0){2}{\vector(0,-1){10}}
            \put(30,30){\vector(1,0){20}}
            \put(30,30){\vector(0,1){20}}
                       
            \put(25,136){\makebox(0,0)[b]{$u_{D}$}}
            \put(50,33){\makebox(0,0)[b]{$x_1$}}
            \put(33,50){\makebox(0,0)[l]{$x_2$}}
            \put(25,0){\makebox(0,0)[t]{1}}
            \put(70,0){\makebox(0,0)[t]{9}}
            \put(142,70){\makebox(0,0)[l]{10}}
            \put(20,4){\line(1,0){100}}
            \put(140,20){\line(0,1){100}}
            \put(20,2){\line(0,1){4}}
            \put(30,2){\line(0,1){4}}
            \put(30,118){\line(0,1){4}}
            \put(120,2){\line(0,1){4}}
            \put(138,20){\line(1,0){4}}
            \put(138,120){\line(1,0){4}}
            
            {\linethickness{2pt} \put(20,120){\line(1,0){10}}}

        \end{picture}
\caption{Geometry of the elastoplastic problem with Drucker-Prager yield criterion.}
\label{fig_scheme2}
\end{figure}

The geometry in $x_1-x_2$ plane is depicted in Figure  \ref{fig_scheme2}. A geometry for the corresponding 3D problem arises from the extrusion in $x_3$ direction of length one. On the left, right and bottom sides of the depicted domain, the zeroth normal displacements are prescribed, i.e., $\mbf u\cdot\mbf n=0$ where $\mbf n$ is a unit normal vector to the boundary. The strip-footing of the length one is considered on the top of the domain. The loading is controlled by the nonhomogeneous Dirichlet boundary condition $u_D$ in the direction $x_2$. The values $u_D$ varies from 0 to 1 using a suitable adaptive strategy described below. Volume and traction forces are not prescribed. Material parameters are set as follows: $E = 1e7\,$ (Young's modulus), $\nu = 0.48$ (Poisson's  ratio), $c_0=450$ (cohesion), and $\phi=\pi/9$ (friction angle). The parameters $c$ and $\eta$ introduced in Section \ref{subsec_DP} are computed by the following formulas \cite[Chapter 6]{NPO08}: 
\begin{eqnarray*}
&&\eta=\frac{6\sin\phi}{\sqrt{3}(3+\sin\phi)}, \quad c=c_0\frac{6\cos\phi}{\sqrt{3}(3+\sin\phi)} \quad \qquad \mbox{for the 3D problem}, \\
&&\eta=\frac{3\tan\phi}{\sqrt{9+12\tan^2\phi}}, \quad c=c_0\frac{3}{\sqrt{9+12\tan^2\phi}} \qquad \mbox{for the plane strain problem. }
\end{eqnarray*}

We start with the constant increment $\triangle u_D=0.001$. For the solution $\mbf u_k$, we compute the corresponding (average) pressure $\hat p_k$ supported by the footing. It is well known that values of $\hat p_k$ are bounded from above by an unknown limit value. So if we observe that the increment $\triangle \hat p_k$ is sufficiently small then we multiply the increment $\triangle u_D$ by factor two to achieve the prescribed maximal displacement faster. 

\begin{figure}
\center
\includegraphics[width=0.47\textwidth]{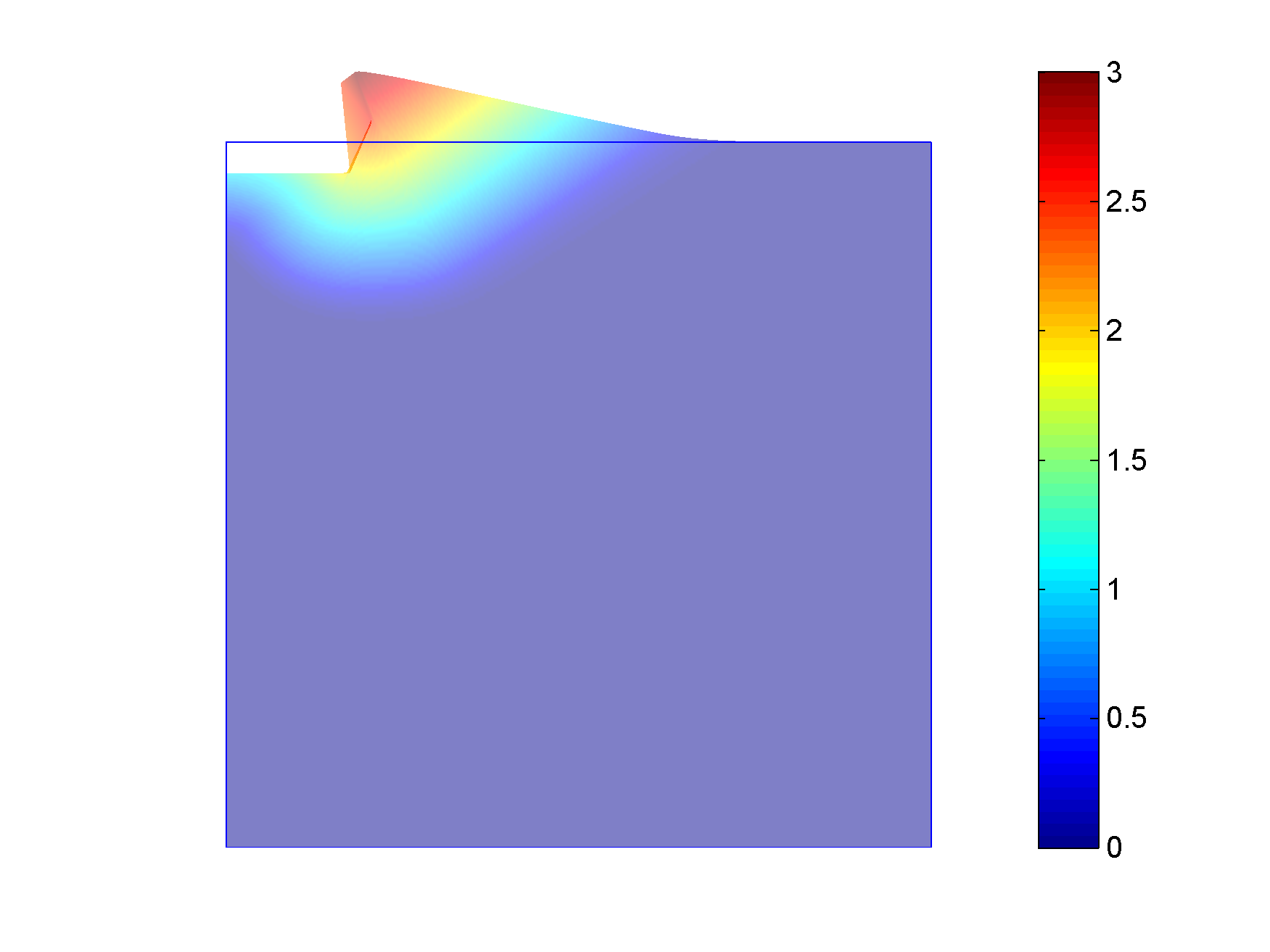}
\includegraphics[width=0.47\textwidth]{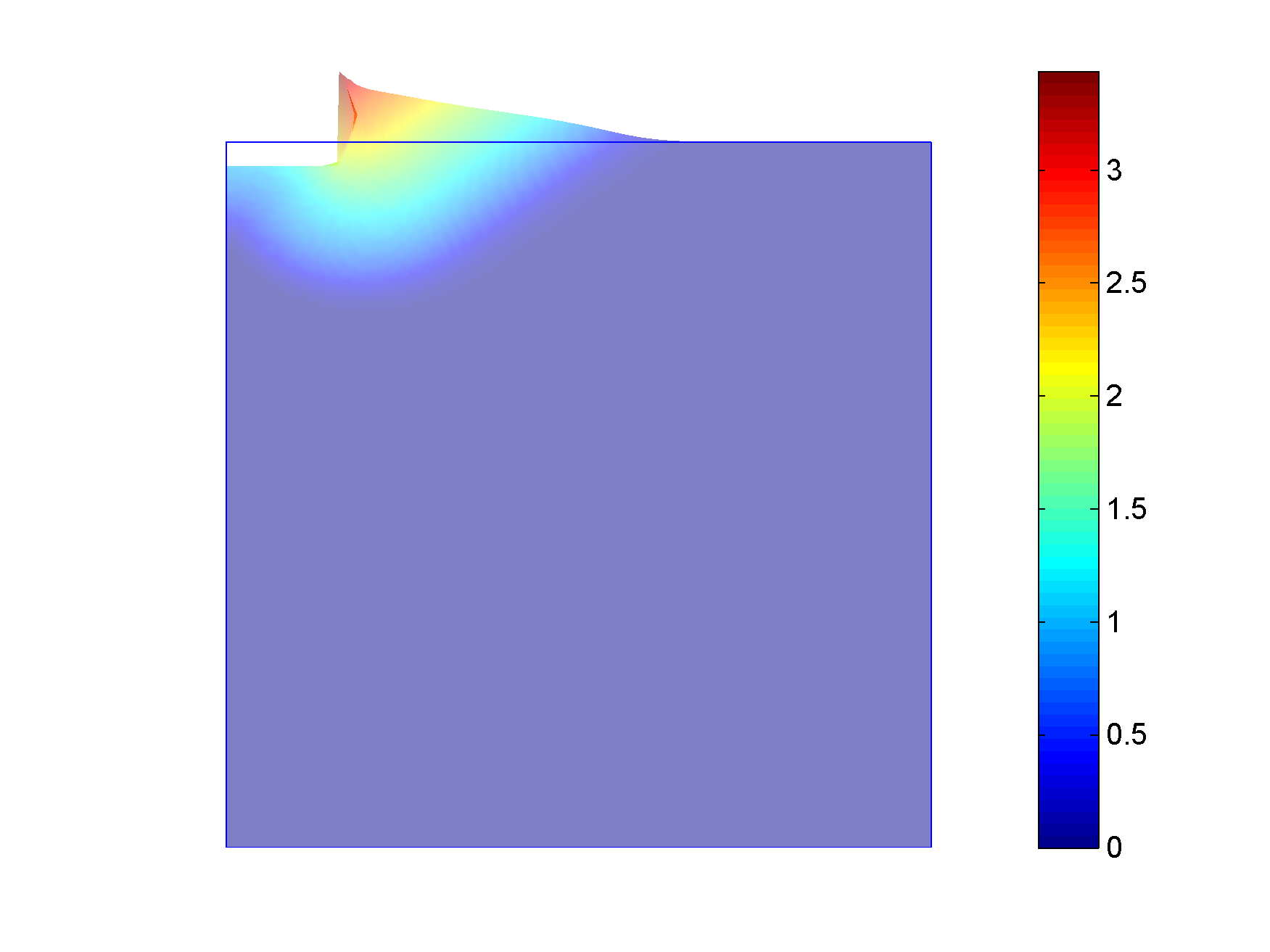}
\caption{Total displacement fields with deform shapes for P1 (left) and P2 (right) elements. The deform shapes correspond to $\mbf u/\max{\mbf u}$.}
\label{fig_DP_deform_shape}
\includegraphics[width=0.45\textwidth]{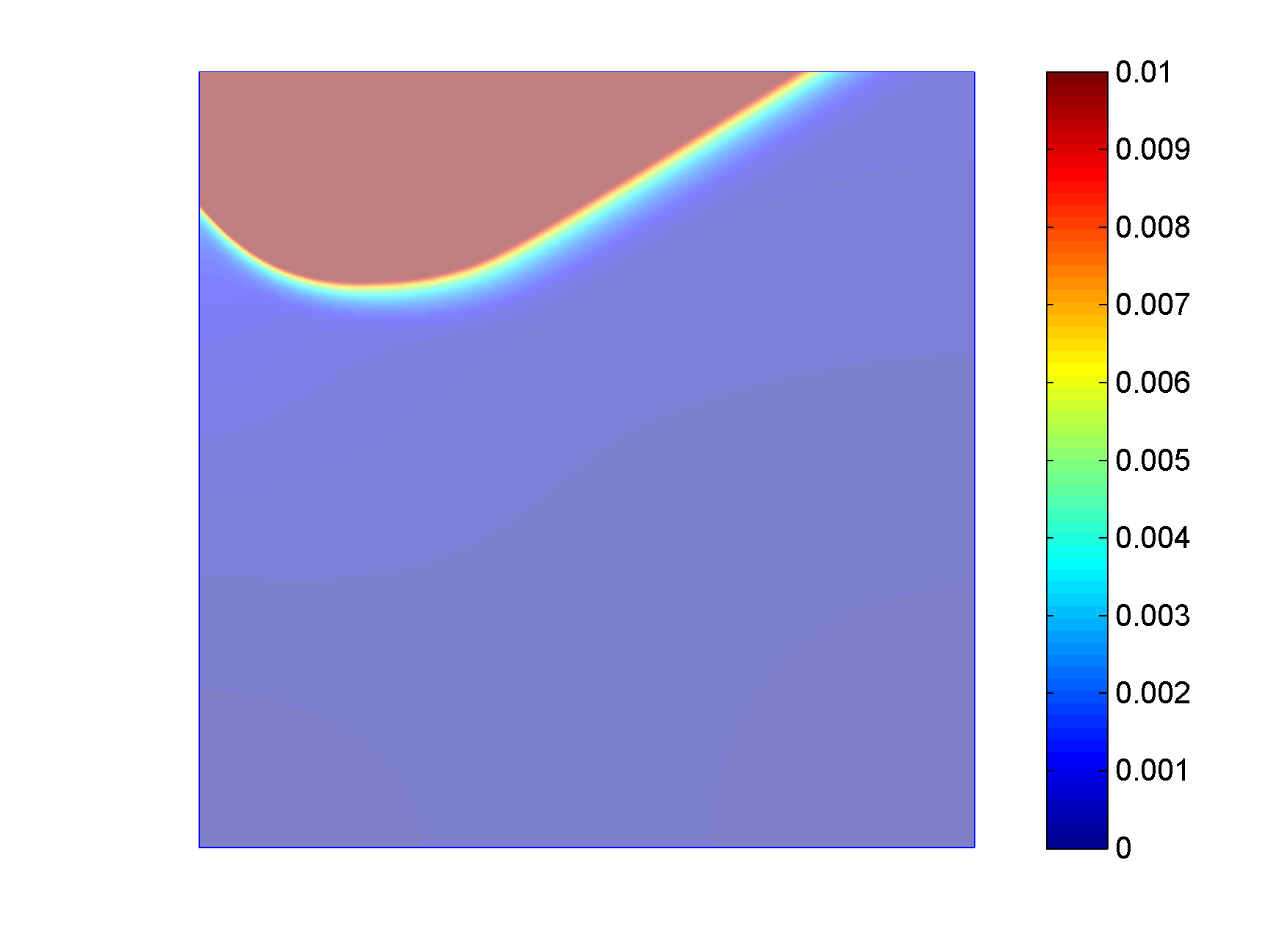}
\includegraphics[width=0.45\textwidth]{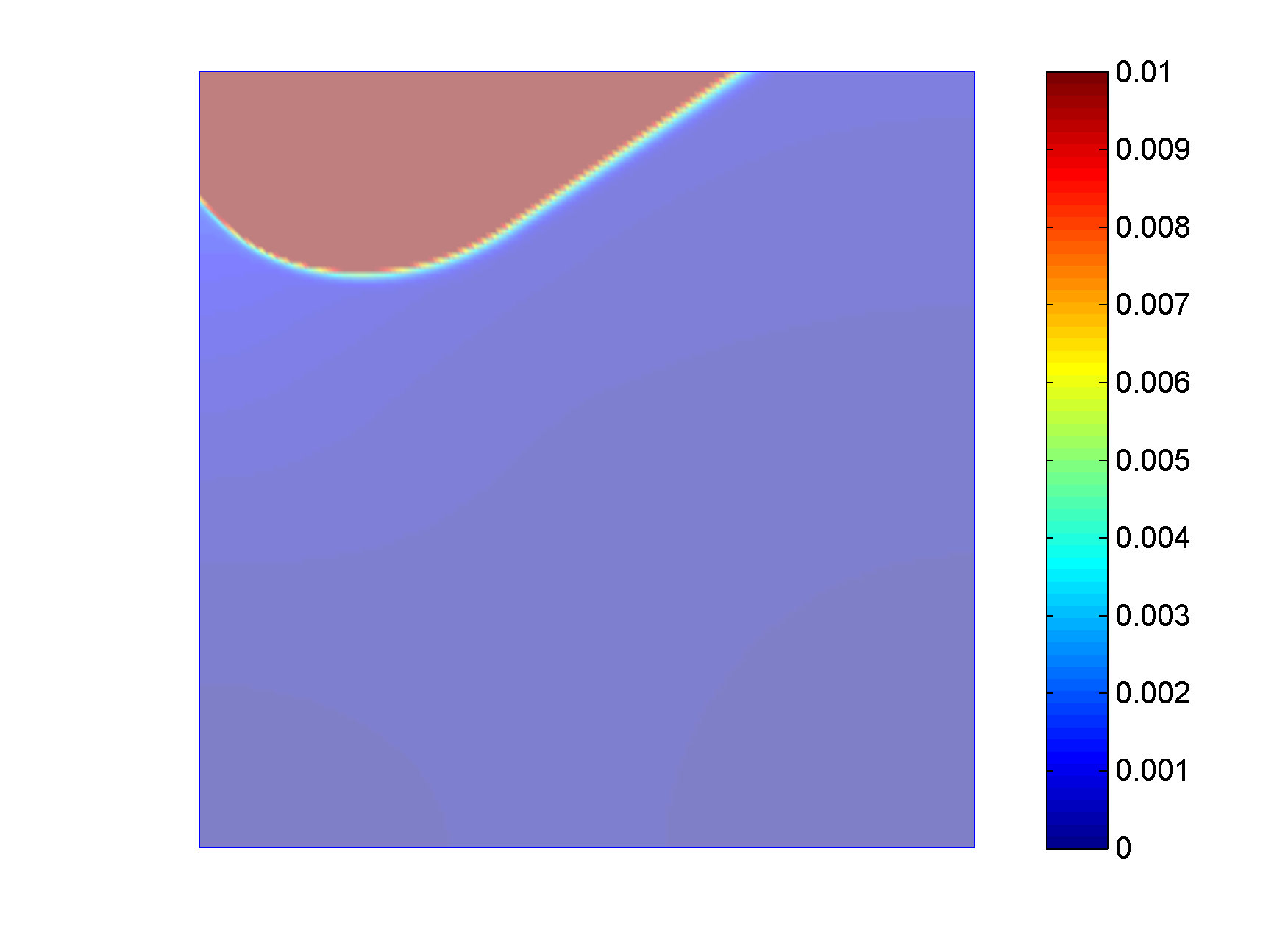}
\caption{Total displacement fields for P1 (left) and P2 (right) elements. Values greater than $0.01$ are replaced with $0.01$.}
\label{fig_DP_displacement}
\includegraphics[width=0.45\textwidth]{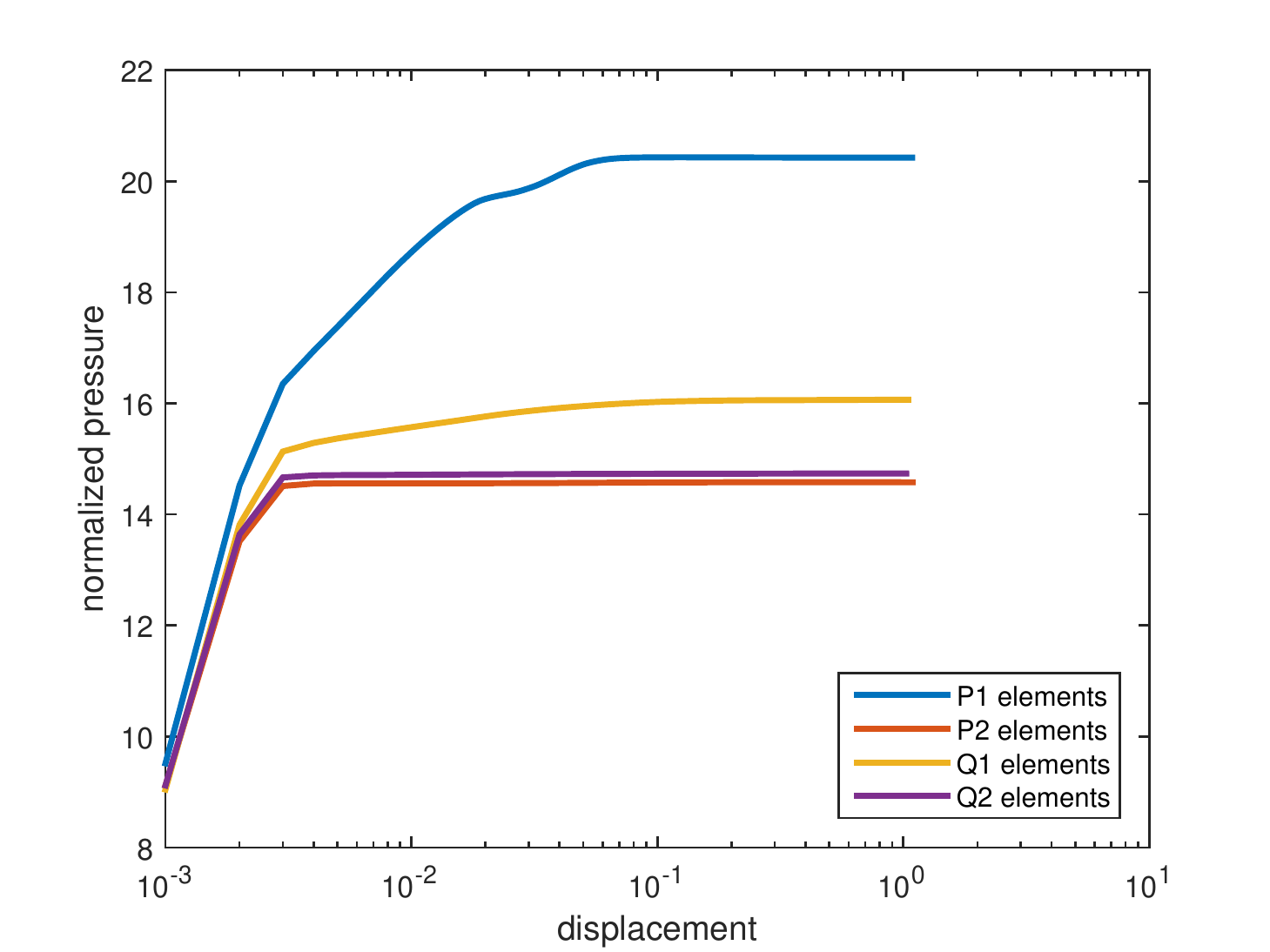}
\includegraphics[width=0.45\textwidth]{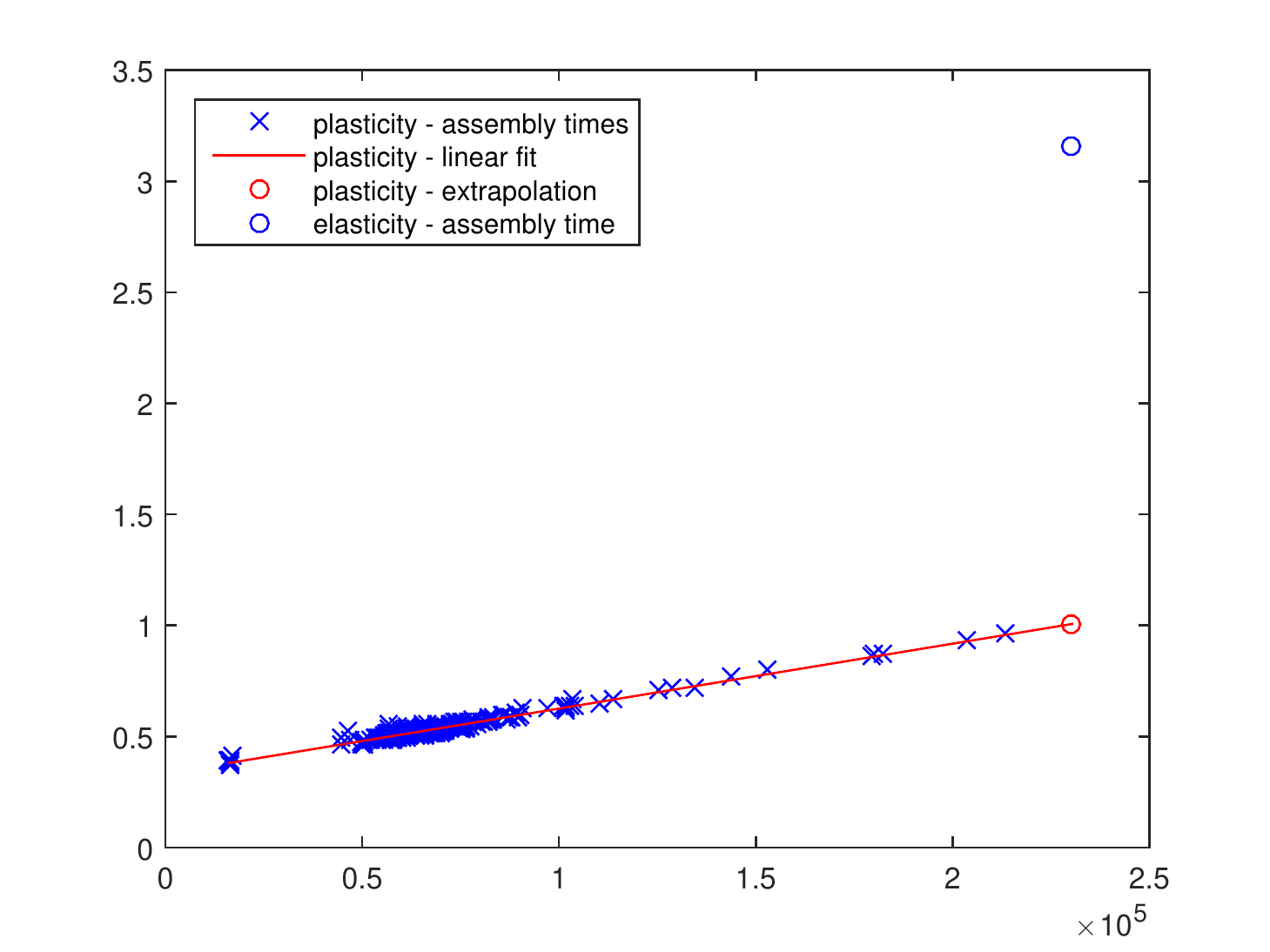}
\caption{Load paths for P1, P2, Q1, Q2 elements (left) and assembly times of tangential stiffness matrix versus number of plastic integration points (right).}
\label{fig_DP_load_path}
\end{figure}

Since the expected results are strongly dependent on the mesh density and on chosen element types, we present the results for 2D problem to achieve finer meshes easily. In particular, we use regular meshes divided the domain into 320$\times$320 squares for P1 and Q1 elements, respectively 160$\times$160 squares for P2 and Q2 elements to have a similar number of unknowns. Figures \ref{fig_DP_deform_shape} and \ref{fig_DP_displacement} compare plastic collapse for P1 and P2 elements. In Figure \ref{fig_DP_deform_shape}, we see total displacement fields with deform shapes. To visualize expected slip surfaces, values of displacements greater than $0.01$ are replaced with $0.01$, see Figure \ref{fig_DP_displacement}. The strong dependence on element types is illustrated in Figure \ref{fig_DP_load_path} (left). We see that the normalized pressures $\hat p_k/c_0$ are significantly overestimated for P1 and Q1 elements. The results for P2 and Q2 elements are in accordance with \cite[Chapter 8]{NPO08}. The assembly times of the tangent stiffness matrices are also illustrated in Figure \ref{fig_DP_load_path} (right). We have a similar observation as for the von Mises yield criterion.

\section{Comment on the technique of P. Byczanski}
Although the basic split of the tangential stiffness matrix \eqref{K-split} can be written as 
\begin{equation} \label{K-split_Petr}
\mbf K_{tangent}=\mbf B^\top\mbf D_{tangent} \mbf B,
\end{equation}

$$$$ this simpler form is not convenient in MATLAB since the difference $\mbf D_{tangent}-\mbf D_{elast}$ of {\it sparse} matrices $\mbf D_{tangent}$ and $\mbf D_{elast}$ can be much sparser than $\mbf D_{tangent}$. This occurs when most of integration points remains in the elastic phase. Therefore, for problems with smaller plastic regions, the assembly of the tangential stiffness matrix can be significantly faster than for problems with larger plastic regions. Let us note that \eqref{K-split_Petr}  was originally used in unpublished codes from P. Byczanski (UGN Ostrava) applied in numerical examples from \cite{Sy12,HRS16b, HRS16a, SCKKZB16, SCL17}, but this idea has not been neither emphasized nor described.

\section{Conclusion and future plans}

The paper is focused on an efficient and flexible implementation of various elastoplastic problems. We have mainly proposed the innovative MATLAB assembly of elastoplastic FEM matrices based on the split (\ref{K-split}). Our techniques are explained and implemented in the vectorized code available for download \cite{software}. Time performance of FEM assembly is comparable with other techniques for purely elastic stiffness matrices. Additional effort to build the tangential stiffness matrices in each Newton iteration and each time step of  elastoplastic problems does not exceed the cost for the elastic stiffness matrix. The smaller is the number of the plastic integrations points, the faster is the assembly. The code is flexible due to the fact that one can choose several types of (Lagrange) finite elements. It can be further extended for various elastoplastic models by changing the function \verb+constitutive_problem+.

\bigskip\noindent
{\bf Acknowledgement.} We would like to thank to our former colleague Petr Byczanski (Ostrava) whose ideas inspired this work. 
The first two authors acknowledge a support from The Ministry of Education, Youth and Sports of the Czech Republic - project LO1404 - Sustainable development of CENET "CZ.1.05/2.1.00/19.0389: Research Infrastructure Development of the CENET" (the first author); project ``IT4Innovations excellence in science - LQ1602" from the National Programme of Sustainability II  (the second author). The third author acknowledges a support by the Czech Science Foundation (GA\v{C}R)  through projects No. 16-34894L, 17-04301S, 18-03834S.

\section{Appendix}

This appendix completes Section \ref{sec_solution} with particular examples of finite elements and related numerical quadratures. We consider P1, P2 tetrahedral elements and Q1, Q2 hexahedral elements for displacement approximation in our implementation.  Higher order P- and Q-type elements can be implemented analogously \cite{So}. We recapitulate these elements for 3D case only, see Figure \ref{fig:3D_elements} for illustration.

\subsection{P1 and P2 tetrahedral elements} 

The reference P1 element is defined on a reference tetrahedron with 4 nodes
$$\hat N_1=[0,0,0], \quad \hat N_2=[1,0,0], \quad \hat N_3=[0,1,0], \quad \hat N_4=[0,0,1]$$
and 4 corresponding linear basis functions (therefore $n_p=4$) are 
$$\hat\Phi_1(\mbf\xi)=1-\xi_1-\xi_2-\xi_3, \quad  \hat\Phi_2(\mbf\xi)=\xi_1, \quad \hat\Phi_3(\mbf\xi)=\xi_2, \quad \hat\Phi_4(\mbf\xi)=\xi_3.$$
Notice that strain fields for $P1$ elements are constant on elements. Therefore, it sufficies to consider 1-point Gauss quadrature, i.e., $n_q=1$, $\hat A_1=[1/4, 1/4, 1/4]$ and $\omega_1=1/6$. We see that the weight coefficient coincides with the volume of the reference element.

The reference P2 element is defined on the same reference tetrahedron above with the nodes $\hat N_1$, $\hat N_2$, $\hat N_3$, $\hat N_4$ and also utilizes 6 edges midpoints
  $$\hat N_5=[1/2,0,0], \quad \hat N_6=[1/2,1/2,0], \quad \hat N_7=[0,1/2,0], \quad \hat N_8=[1/2,0,1/2], \quad \hat N_9=[0,1/2,1/2], \quad \hat N_{10}=[0,0,1/2].$$
Let $\xi_0:=\xi_0(\mbf\xi)=1-\xi_1-\xi_2-\xi_3$. Then the quartet $(\xi_0,\xi_1,\xi_2,\xi_3)$ defines the barycentric coordinates and one can write the quadratic basis functions as follows:
$$
\begin{array}{c}
\hat\Phi_1(\mbf\xi)=\xi_0(2\xi_0-1),\;\; \hat\Phi_2(\mbf\xi)=\xi_1(2\xi_1-1),\;\; \hat\Phi_3(\mbf\xi)=\xi_2(2\xi_2-1),\;\; \hat\Phi_4(\mbf\xi)=\xi_3(2\xi_3-1),\\[2mm]
\hat\Phi_5(\mbf\xi)=4\xi_0\xi_1,\;\;\hat\Phi_6(\mbf\xi)=4\xi_1\xi_2,\;\;\hat\Phi_7(\mbf\xi)=4\xi_0\xi_2,\;\;\hat\Phi_8(\mbf\xi)=4\xi_1\xi_3,\;\;\hat\Phi_9(\mbf\xi)=4\xi_2\xi_3,\;\;\hat\Phi_{10}(\mbf\xi)=4\xi_0\xi_3.
\end{array}
$$
For P2 elements, we use 11-point numerical quadrature which is exact to order 4. The coordinates of the quadrature points and their weights are following \cite{Yu84}:
$$
\begin{array}{ll}
\hat A_1\;=[0.250000000000000,\;  0.250000000000000,\; 0.250000000000000],& \omega_1=-0.013155555555555,\\
\hat A_2\;=[0.071428571428571,\;  0.071428571428571,\; 0.071428571428571],& \omega_2\;=\;0.007622222222222,\\
\hat A_3\;=[0.785714285714286,\;  0.071428571428571,\; 0.071428571428571],& \omega_3\;=\;0.007622222222222,\\
\hat A_4\;=[0.071428571428571,\;  0.785714285714286,\; 0.071428571428571],& \omega_4\;=\;0.007622222222222,\\
\hat A_5\;=[0.071428571428571,\;  0.071428571428571,\; 0.785714285714286],& \omega_5\;=\;0.007622222222222,\\
\hat A_6\;=[0.399403576166799,\;  0.100596423833201,\; 0.100596423833201],& \omega_6\;=\;0.024888888888888,\\
\hat A_7\;=[0.100596423833201,\;  0.399403576166799,\;  0.100596423833201],& \omega_7\;=\;0.024888888888888,\\
\hat A_8\;=[0.100596423833201,\; 0.100596423833201, \; 0.399403576166799 ],& \omega_8\;=\;0.024888888888888,\\
\hat A_9\;=[0.399403576166799,\;  0.399403576166799,\; 0.100596423833201],& \omega_9\;=\;0.024888888888888,\\
\hat A_{10}=[0.399403576166799,\;  0.100596423833201,\; 0.399403576166799],& \omega_{10}=\;0.024888888888888,\\
\hat A_{11}=[0.100596423833201,\; 0.399403576166799,\;  0.399403576166799],& \omega_{11}=\;0.024888888888888.
\end{array}
$$

\subsection{Q1 and Q2 hexahedral elements}
The reference Q1 elements  is defined on a hexahedron with 8 nodes
\begin{eqnarray*}
&&\hat N_1=[-1,-1,-1], \quad \hat N_2=[1,-1,-1], \quad \hat N_3=[1,1,-1], \quad \hat N_4=[-1,1,-1], \\
 &&\hat N_5=[-1,-1,1], \quad \hat N_6=[1,-1,1], \quad \hat N_7=[1,1,1], \quad \hat N_8=[-1,1,1]
\end{eqnarray*}
and the corresponding linear basis functions  ($n_p=8$) are
\begin{eqnarray*}
\hat\Phi_1(\mbf\xi)=\frac{1}{8}(1-\xi_1)(1-\xi_2)(1-\xi_3),& \hat\Phi_5(\mbf\xi)=\frac{1}{8}(1-\xi_1)(1-\xi_2)(1+\xi_3),\\
\hat\Phi_2(\mbf\xi)=\frac{1}{8}(1+\xi_1)(1-\xi_2)(1-\xi_3),& \hat\Phi_6(\mbf\xi)=\frac{1}{8}(1+\xi_1)(1-\xi_2)(1+\xi_3),\\
\hat\Phi_3(\mbf\xi)=\frac{1}{8}(1+\xi_1)(1+\xi_2)(1-\xi_3),& \hat\Phi_7(\mbf\xi)=\frac{1}{8}(1+\xi_1)(1+\xi_2)(1+\xi_3),\\
\hat\Phi_4(\mbf\xi)=\frac{1}{8}(1-\xi_1)(1+\xi_2)(1-\xi_3),& \hat\Phi_8(\mbf\xi)=\frac{1}{8}(1-\xi_1)(1+\xi_2)(1+\xi_3).
\end{eqnarray*}
We use $2\times2\times2$ Gauss quadrature derived from the 1D case where the quadrature points are located at $-1/\sqrt{3}$ and $1/\sqrt{3}$, and the corresponding weights are equal to one  (see \cite{B06}).

The reference Q2 element is defined on the same reference hexahedron above and also utilizes 12 edges midpoints 
\begin{eqnarray*}
\hat N_9=[0,-1,-1], \hat N_{10}=[1,0,-1],  \quad \hat N_{11}=[0,1,-1],  \quad \hat N_{12}=[-1,0,-1],  \quad \hat N_{13}=[0,-1,1],  \quad \hat N_{14}=[1,0,1], \\
\hat N_{15}=[0,1,1],  \quad \hat N_{16}=[-1,0,1],  \quad \hat N_{17}=[-1,-1,0],  \quad \hat N_{18}=[1,-1,0],  \quad \hat N_{19}=[1,1,0],  \quad \hat N_{20}=[-1,1,0].
\end{eqnarray*}
The corresponding quadrature basis functions are defined as follows:
$$
\begin{array}{rclcrcl}
\hat\Phi_1(\mbf\xi) & = & \frac{1}{8}(1-\xi_1)(1-\xi_2)(1-\xi_3)(-2-\xi_1-\xi_2-\xi_3),&& \hat\Phi_5(\mbf\xi) & = & \frac{1}{8}(1-\xi_1)(1-\xi_2)(1+\xi_3)(-2-\xi_1-\xi_2+\xi_3),\\
\hat\Phi_2(\mbf\xi) & = & \frac{1}{8}(1+\xi_1)(1-\xi_2)(1-\xi_3)(-2+\xi_1-\xi_2-\xi_3),&& \hat\Phi_6(\mbf\xi) & = &\frac{1}{8}(1+\xi_1)(1-\xi_2)(1+\xi_3)(-2+\xi_1-\xi_2+\xi_3),\\
\hat\Phi_3(\mbf\xi) & = & \frac{1}{8}(1+\xi_1)(1+\xi_2)(1-\xi_3)(-2+\xi_1+\xi_2-\xi_3),&& \hat\Phi_7(\mbf\xi) & = & \frac{1}{8}(1+\xi_1)(1+\xi_2)(1+\xi_3)(-2+\xi_1+\xi_2+\xi_3),\\
\hat\Phi_4(\mbf\xi) & = & \frac{1}{8}(1-\xi_1)(1+\xi_2)(1-\xi_3)(-2-\xi_1+\xi_2-\xi_3),&& \hat\Phi_8(\mbf\xi) & = & \frac{1}{8}(1-\xi_1)(1+\xi_2)(1+\xi_3)(-2-\xi_1+\xi_2+\xi_3),\\
\hat\Phi_9(\mbf\xi) & = & \frac{1}{4}(1-\xi_1^2)(1-\xi_2)(1-\xi_3),   && \hat\Phi_{10}(\mbf\xi) & = & \frac{1}{4}(1+\xi_1)(1-\xi_2^2)(1-\xi_3),\\
\hat\Phi_{11}(\mbf\xi) & = & \frac{1}{4}(1-\xi_1^2)(1+\xi_2)(1-\xi_3),&& \hat\Phi_{12}(\mbf\xi) & = & \frac{1}{4}(1-\xi_1)(1-\xi_2^2)(1-\xi_3),\\
\hat\Phi_{13}(\mbf\xi) & = & \frac{1}{4}(1-\xi_1^2)(1-\xi_2)(1+\xi_3),&&  \hat\Phi_{14}(\mbf\xi) & = & \frac{1}{4}(1+\xi_1)(1-\xi_2^2)(1+\xi_3),\\
\hat\Phi_{15}(\mbf\xi) & = & \frac{1}{4}(1-\xi_1^2)(1+\xi_2)(1+\xi_3),&& \hat\Phi_{16}(\mbf\xi) & = & \frac{1}{4}(1-\xi_1)(1-\xi_2^2)(1+\xi_3),\\
\hat\Phi_{17}(\mbf\xi)& = & \frac{1}{4}(1-\xi_1)(1-\xi_2)(1-\xi_3^2),&& \hat\Phi_{18}(\mbf\xi) & = & \frac{1}{4}(1+\xi_1)(1-\xi_2)(1-\xi_3^2),\\
\hat\Phi_{19}(\mbf\xi) & = & \frac{1}{4}(1+\xi_1)(1+\xi_2)(1-\xi_3^2),&& \hat\Phi_{20}(\mbf\xi) & = & \frac{1}{4}(1-\xi_1)(1+\xi_2)(1-\xi_3^2).
\end{array}
$$
We use $3\times3\times3$ Gauss quadrature derived from the 1D case where the quadrature points are located at $-\sqrt{3/5}$, $0$, and  $\sqrt{3/5}$, and the corresponding weights are equal to 5/9, 8/9, and 5/9, respectively (see \cite{B06}).




\section*{References}
\bibliographystyle{plain}
\bibliography{Ce_Sy_Va_biblio}

\end{document}